\documentclass[draft,12pt]{article}
\usepackage{amssymb}
\usepackage[final]{graphicx}
\usepackage{epsfig}
\newtheorem{defn}{Definition}[section]
\newtheorem{lemma}[defn]{Lemma}

\newtheorem{theorem}[defn]{Theorem}
\newtheorem{definition}[defn]{Definition}

\newtheorem{prop}[defn]{Proposition}
\newtheorem{corollary}[defn]{Corollary}

\newcommand{\C}{\mathbb{C}}
\newcommand{\Z}{\mathbb{Z}}
\newcommand{\R}{\mathbb{R}}
\newcommand{\Q}{\mathbb{Q}}
\newcommand{\Cn}{\mathbb{C}^{n}}
\newcommand{\aand}{\textrm{ and }}
\newcommand{\bdm}{\begin{displaymath}}
\newcommand{\edm}{\end{displaymath}}
\newcommand{\beqn}{\begin{equation}}
\newcommand{\eeqn}{\end{equation}}
\newcommand{\partialf}[1]{\frac{\partial f(#1)}{\partial z_{j}}}
\newcommand{\partialfi}[2]{\frac{\partial f_{#1}(#2)}{\partial z_{j}}}
\newcommand{\Bhh}{B_{|hh|}(\tilde{a})}
\newcommand{\bbb}{\textrm{   }}
\newcommand{\bb}{\textrm{  }}
\newcommand{\bs}{\textrm{ }}
\newcommand{\bbbbb}{\bbb \bbb \bbb \bbb \bbb \bbb \bbb \bbb \bbb \bbb}
\newcommand{\BBb}{\bbbbb \bbbbb} 
\newcommand{\BbB}{\bbb \bbb \bbb \bbb \bbb \bbb \bbb \bbb \bbb \bbb \bbb \bbb \bbb \bbb \bbb}
\newcommand{\eol}{$\setminus$}
\newcommand{\eolm}{\setminus$}
\newcommand{\Aff}{\textrm{Aff}}
\newcommand{\llg}{\, \textrm{log}}
\newcommand{\Lsigv}[1]{L_{\sigma_{#1}}(v)}
\newcommand{\Lsig}[1]{L_{\sigma_{#1}}}
\newcommand{\Jsig}[1]{J_{\sigma_{#1}}} 
\newcommand{\J}{\mathcal{J}}
\newcommand{\Rm}{\mathbf{R}}
\newcommand{\Rmo}{\overline{\mathbf{R}}}  
\newcommand{\Rmu}{\underline{\mathbf{R}}}  
\newcommand{\Cm}{\mathbf{C}}
\newcommand{\Rmb}{\mathbf{\underline{R}_{\beta}}} 
\newcommand{\Rb}{\mathbf{R_{\beta}}}  
\newcommand{\Mm}{\mathbf{M}}
\newcommand{\Mhm}{\mathbf{M_{h}}}
 
\newcommand{\Mhmu}{\mathbf{\underline{M}_{h}}} 
\newcommand{\Lm}{\mathbf{L}}
\newcommand{\MLm}{\mathbf{ML}}
\newcommand{\Sm}{\mathbf{S}}
\newcommand{\Smu}{\underline{\mathbf{S}}} 
\newcommand{\Fm}{\mathbf{F}} 
\newcommand{\Fmu}{\underline{\mathbf{F}}} 
\newcommand{\MLu}{\underline{\mathbf{ML}}}
\newcommand{\im}{\textrm{im}} 
\renewcommand{\ker}{\textrm{ker}} 

\newcommand{\bqr}{\begin{eqnarray*}} 
\newcommand{\eqr}{\end{eqnarray*}}
\newcommand{\hgamma}{\hat{\gamma}} 
\newcommand{\hdelta}{\hat{\delta}}
\begin{document}
\title{Proving A Manifold To Be Hyperbolic Once It Has\\
Been Approximated To Be So}
\author{Harriet Handel Moser}
\maketitle
\pagenumbering{arabic}
\section{Introduction} \label{intro}
This paper presents the major result of my doctoral dissertation written at Columbia University~\cite{HMdis}, with Walter Neumann as my thesis adviser.  Known uses of the method developed, which allows one to conclusively prove that a 3-manifold has a hyperbolic structure, include some of David Gabai's~\cite{Gab1,Gab2} recent work and a paper by Chris Leininger~\cite{Lein}.  Since the determination that $M$ is complete hyperbolic is dependent on there being a solution to a set of equations, we shall first review the development of these equations.  Every orientable complete hyperbolic manifold of finite volume is obtained from an ideally triangulated one by Dehn surgery on some of its cusps.  This fact is documented by Neumann~\cite{N-Z}, based on a Thurston preprint~\cite{Thur2}, so we first examine $N$, a non-compact 3-manifold that is the interior of a compact one whose boundary consists of $k$ tori.  $N$ can be realized as a gluing of $n$ tetrahedra, $\sigma_{1},\ldots,\sigma_{n}$, having $k$ vertices after gluing, with a \emph{conic neighborhood} of each vertex removed~\cite{B-P}.  A \emph{conic neighborhood of the vertex, $v,$} is described as follows.  Let $v$ be a vertex and $\sigma_{j}$ a tetrahedron that $v$ belongs to.  Take the second barycentric subdivision of the edges of $\sigma_{j}$ containing $v$ and let $w_{1},\, w_{2} \aand w_{3}$ be the closest vertices to $v$ for these edges with respect to this subdivision.  See Figure~\ref{tet1}.
\begin{figure}
\begin{picture}(160,120)(-170,0)
\put(-19,0){\line(6,0){114}}
\put(95,0){\line(-1,2){57}}
\put(-19,0){\line(1,2){57}}
\put(-19,0){\line(2,1){57}}
\put(95,0){\line(-2,1){57}}
\put(9.5,57){\line(1,0){57}}
\put(66.5,57){\line(-2,1){28.5}}
\put(9.5,57){\line(2,1){28.5}}
\put(38,114){\line(0,-1){39}}
\put(38,27.5){\line(0,1){25}}
\put(-9,54){$w_{1}$}
\put(73,54){$w_{2}$}
\put(39,76){$w_{3}$}
\put(38,120){$v$}
\put(106,65){\vector(-1,0){67}}
\put(107,65){$L_{\sigma_{j}}(v)$}
\end{picture}
\caption{The Tetrahedron $\sigma_{j}$}  \label{tet1}
\end{figure}
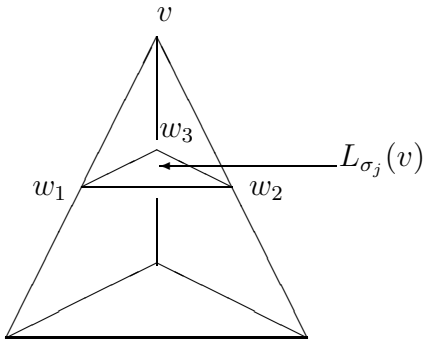
\begin{definition}
\begin{itemize}
\item $\mathbf{L_{\sigma_{j}}(v)} = \mathrm{triangle \ having \ vertices \ } w_{1},\, w_{2} \aand w_{3} \ \mathrm{as}$ 

$\mathrm{above \ with \ respect \ to \ }v \aand \sigma_{j}$ 
\item $\mathbf{L(v)} = \bigcup_{\stackrel{v \mathrm{\, vertex \,of \,}\sigma_{j}}{\scriptscriptstyle 1\leq j\leq n}} L_{\sigma_{j}}(v)$

$L(v)$ is called the $\mathbf{link}$ of $v$
\item $\mathbf{U_{\sigma_{j}}(v)} = \mathrm{tetrahedron \ having \ vertices \ } v,\, w_{1},\, w_{2} \aand w_{3}$
\item $\mathbf{conic \ neighborhood \ of \ }\mathbf{v} = \bigcup_{\stackrel{v \mathrm{\, vertex \,of \,}\sigma_{j}}{\scriptscriptstyle 1\leq j\leq n}} U_{\sigma_{j}}(v).$
\end{itemize}
\end{definition}

Every vertex is identified with a cusp of $N,$ and its link is a torus.  These truncated tetrahedra resulting from the removal of the conic sections can now be treated as ideal hyperbolic tetrahedra, so there exists a hyperbolic structure on $N\backslash \textrm{1-skeleton of }N.$  In order for $N$ to have a hyperbolic structure, there must be consistency across the 1-skeleton.  The conditions for this to happen are embodied in the consistency equations and will be described in detail in Section~\ref{ideq}, ``Identifying the Equations". 

Completeness applies to the cusps.  Once a hyperbolic structure is identified, it induces a \emph{similarity structure} (i.e., a $(\C, \textrm{Aff}(\C))$ structure) on each of the $k$ tori, $T_{1},\ldots,T_{k}$.  If the similarity structure of a torus identified with a cusp is Euclidean, $N$ will be complete at that cusp~\cite{B-P}.  This occurs when the image of the holonomy of the similarity structure for the torus consists entirely of translations, or equivalently, has at least one non-trivial translation~\cite{B-P}.  A \emph{holonomy} of a similarity structure for a torus, $T,$ is a map $\theta$ such that \mbox{$\theta:\pi_{1}(T)\to \textrm{Aff}(\C)$~\cite{B-P}.} The conditions for the image of $\theta$ to consist entirely of translations are presented by the completeness equations which will also  be discussed in Section~\ref{ideq}, ``Identifying The Equations".

Once we establish the conditions for cusps of $N$ to be complete, we turn our attention to the manifold $M$, obtained from $N$ by Dehn surgery on some of the cusps.  Assume $h$ cusps remain unsurgered, so there are $k-h$ surgered cusps.  $M$  must satisfy the consistency equations; however, there are now only $h$ cusps that must be shown to be complete, so we only need the completeness equations referring to these $h$ cusps.  The remaining $k-h$ surgered cusps must result from Dehn surgery with co-prime coefficients $(p_{i},q_{i}) \textrm{ for }1\leq i\leq k-h$ where $(p_{i},q_{i})$ and the holonomy of the similarity structure of $T_{i}$ are joined in one equation~\cite{B-P}.

Once the equations needed to prove a manifold complete hyperbolic are identified, we set up the machinery to test whether a solution exists in Section~\ref{sol}, ``How to Test for a Solution."  The method described there concludes the proof of the following theorem, which is our main result.

\begin{theorem} \label{th}
Let $M$ be a manifold and assume there are $n$ tetrahedra in the triangulation of $M$ according to SnapPea~\cite{Snapp}.  There are $n$ equations, $\{ f_{i}(z)=0 \, | \, f_{i}:\Cn \to \C \}$ for $1 \leq i \leq n,$ whose simultaneous solution will guarantee that $M$ is complete hyperbolic.  If SnapPea finds an approximate geometric solution to these equations, let $a = (a_{1}, \ldots, a_{n})$ be an approximate geometric solution generated by SNAP~\cite{Snap} on the SnapPea manifold file for $M.$  Let $b_{i} = f_{i}(a)$ for $1 \leq i \leq n$ and $f:\Cn \to \Cn$ with $f(z) = (f_{1}(z), \ldots, f_{n}(z)),$ so $f(a) = b = (b_{1}, \ldots, b_{n})$.  Then there is $L>0$ such that there is a genuine solution to the equations, making $M$ complete hyperbolic when the following inequality is true:
\bdm
\textrm{  }|b| \leq \frac{1}{2L|f'(a)^{-1}|^{2}}.
\edm
\end{theorem}

We devote the final section to examples.  Every manifold in the cusped census of SnapPea has been examined and the results are reported in Section~\ref{ex}, ``Examples."  However, for detailed discussion, three examples are presented.  There are simple ones, such as the figure $8$ knot complement and Dehn surgery on the Whitehead link complement.  There is also a complicated link complement with 4 cusps and 32 tetrahedra.  In uncomplicated cases, it is sometimes possible to show that a knot or link complement has a complete hyperbolic structure using means other than the SnapPea approximation.  Thurston has proven that the figure 8 knot complement has a complete hyperbolic structure, and shown when a $(p,q)$ Dehn filling has the same property~\cite{Thur}.  Neumann and Reid have done the same for Dehn fillings of the Whitehead link~\cite{N-R}.  However, when it comes to complicated knots and links, until now, it may have been impossible to definitively determine whether this structure exists.  For several years Leininger had withheld publication of his paper devoted two very large links, one of which is the last example~\cite{Lein}, because he could not prove that their complements have a complete hyperbolic structure.  The paper has now been released using the method presented here.  So far, every manifold that has an approximate solution with respect to a geometric triangulation in SnapPea that has been tested by this method has been verified to have a complete hyperbolic structure.

\section{Identifying the Equations} \label{ideq}
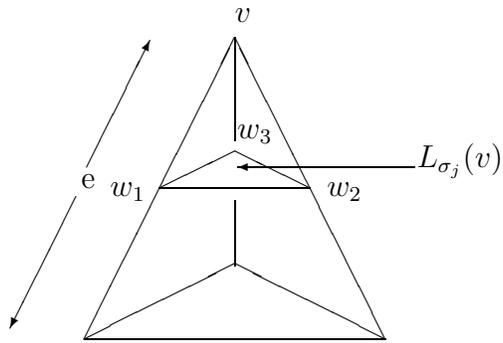
\begin{figure}
\begin{picture}(160,120)(-170,0)
\put(-19,0){\line(6,0){114}}
\put(95,0){\line(-1,2){57}}
\put(-19,0){\line(1,2){57}}
\put(-19,0){\line(2,1){57}}
\put(95,0){\line(-2,1){57}}
\put(9.5,57){\line(1,0){57}}
\put(66.5,57){\line(-2,1){28.5}}
\put(9.5,57){\line(2,1){28.5}}
\put(38,114){\line(0,-1){39}}
\put(38,27.5){\line(0,1){25}}
\put(-9,54){$w_{1}$}
\put(73,54){$w_{2}$}
\put(39,76){$w_{3}$}
\put(38,120){$v$}
\put(106,65){\vector(-1,0){67}}
\put(107,65){$L_{\sigma_{j}}(v)$}
\put(-20,57){e}
\put(-18,65 ){\vector(1,2){24}}
\put(-22,54){\vector(-1,-2){25}}
\end{picture}
\caption{Edge $e$ of the Tetrahedron $\sigma_{j}$}  \label{tet2}
\end{figure}

Let $\sigma_{j}$ be an  ideal hyperbolic tetrahedron as described in Section~\ref{intro}, ``Introduction", and pick an edge $e$ such that $w_{1}\in e$ and prior to truncation, $e$ ended in the vertex $v,$ as in Figure~\ref{tet2}.  Then $L_{\sigma_{j}}(v), $ the triangle with vertices $w_{1},\, w_{2} \aand w_{3}$ naturally has a similarity structure as the triangle in $\C$ with vertices $0,\, 1 \aand z$ (see Figure~\ref{triangle})~\cite{Thur1,N-Z,B-P}, and the dihedral angle at $e$ will be arg$(z).$  Clearly, $z$ must be in $\C_{+},$ the upper half plane in $\C$.  The \emph{modulus of $L_{\sigma_{j}}(v)$ with respect to $w_{1}$} is $z,$ so that the inner angle of the triangle at $w_{1}$ is arg$(z).$  The \emph{modulus of $\sigma_{j}$ at edge $e$} is $z.$  The only other modulii at the other edges of $\sigma_{j}$ will be either $1-(1/z)$ or $1/(1-z),$ so $z$ uniquely describes $\sigma_{j}$ in the upper half plane.  There are six edges with opposite edges having the same modulus~\cite{N-Z,B-P,Ratc}.  See Figure~\ref{tet3}.
\begin{figure}
\begin{picture}(160,95)(-120,-10)
\put(0,0){\line(2,0){152}}
\put(0,0){\line(1,1){76}}
\put(152,0){\line(-1,1){76}}
\put(-3,-3){$0$}
\put(152,-3){$1$}
\put(76,80){$z$}
\put(12,5){$z$}
\put(120,7){$\frac{1}{1-z}$}
\put(62,53){$1-\frac{1}{z}$}
\end{picture}
\caption{The Triangle Similar to $L_{\sigma_{j}}(v)$}  \label{triangle}
\end{figure}
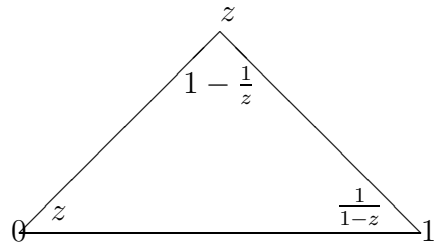
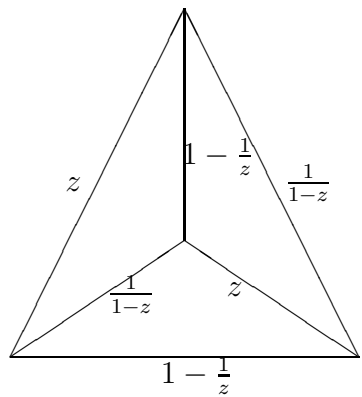
\begin{figure}
\begin{picture}(160,140)(-120,-15)
\put(0,0){\line(4,0){132}}
\put(132,0){\line(-1,2){66}}
\put(0,0){\line(1,2){66}}
\put(0,0){\line(3,2){66}}
\put(132,0){\line(-3,2){66}}
\put(66,44){\line(0,3){88}}
\put(82,23){$z$}
\put(21,63){$z$}
\put(104,63){$\frac{1}{1-z}$}
\put(37,21){$\frac{1}{1-z}$}
\put(65,73){$1-\frac{1}{z}$}
\put(57,-10){$1-\frac{1}{z}$}        
\end{picture}
\caption{Modulii Associated To Edges of the Tetrahedron $\sigma_{j}$}  \label{tet3}
\end{figure} 
\subsection{Consistency Equations}
In order for $N$ to be hyperbolic, if $e$ is an edge of $N,$ the tetrahedra gluing together at $e$ must close up around $e.\,$ That is, the product of all the edge modulii associated with $e$ (different modulus for each tetrahedron $e$ belongs to) must be $e^{2\pi {i}},$ assuring that the sum of the arguments is precisely $2\pi.$  Any of the three distinct edge modulii of a tetrahedron, $\sigma_{j},$ can be expressed as $\pm z_{j}^{r'_{j}}(1-z_{j})^{r''_{j}}$ with $(r'_{j},r''_{j})\in \{(1,0),(-1,1),(0,-1)\},$ so the gluing requirement at edge $e$ is
\begin{displaymath}
\prod_{j=1}^{n} z_{j}^{r'_{j}}(1-z_{j})^{r''_{j}} = \pm 1,
\end{displaymath}
where $r_{j}' = r_{j}'' = 0$ if $\sigma_{j}$ does not contain $e$.  A tetrahedron can have more than one edge glued at $e$ so $r_{j}' \aand r_{j}''$ can take values between $-2 \aand 2$.  The Euler characteristic of $N$ is zero, so it can be shown that $N$ has $n$ edges~\cite{N-Z}.  Thus, the $n$ edge equations can be expressed as 
\beqn
\prod_{j=1}^{n} z_{j}^{r'_{ij}}(1-z_{j})^{r''_{ij}} = \pm 1 \qquad (i=1, \ldots, n).  \label{rval}
\eeqn
They are referred to as the consistency equations.  The existence of a solution is sufficient to make $N$ hyperbolic.  We rewrite them as log equations because they are easier to use  this way and it reflects the fact that the sum of the arguments  of the modulii at each edge is exactly $2\pi$~\cite{N}.
\begin{equation}
\sum_{j=1}^{n}({r'_{ij} \llg(z_{j})} + {r''_{ij} \llg(1-z_{j})}) = c_{i} \pi {i} \qquad c_{i}\in \Z \qquad (i=1, \ldots , n)
\end{equation}
Let $\mathbf{R},\,\mathbf{C} \aand \mathbf{\overline{R}}$ be \label{Rmat} the following matrices.

\begin{displaymath}
\mathbf{R} =
\left( \begin{array}{cccccc}
r'_{11} & \ldots & r'_{1n} & r''_{11} & \ldots & r''_{1n} \\
\vdots & \ddots & \vdots & \vdots & \ddots & \vdots \\
r'_{n1} & \ldots & r'_{nn} & r''_{n1} & \ldots & r''_{nn} 
\end{array} \right) \quad \mathbf{C} = 
\left( \begin{array}{c}
-c_{1} \\ \vdots \\ -c_{n}
\end{array} \right) \quad \mathbf{\overline{R} = (R,C)} 
\end{displaymath}

\begin{prop}
If rank $\mathbf{\overline{R}} = p,$ then the space of solutions to the consistency equations can be defined by exactly $p$ consistency equations.
\end{prop}

\emph{Proof.} Let rank $\mathbf{\overline{R}} = p \leq n,$ so, without loss of generality, we can assume the first $p$ rows of $\mathbf{\overline{R}}$ are linearly independent.  For $s > p,$ there exist $\lambda^{s}_{i} \in \C \textrm{ for }1 \leq i \leq p$ such that

\begin{displaymath}
r'_{sj} = \sum_{i=1}^{p}\lambda^{s}_{i}r'_{ij} \qquad r''_{sj} = \sum_{i=1}^{p}\lambda^{s}_{i}r''_{ij} \qquad c_s = \sum_{i=1}^{p}\lambda^{s}_{i}c_{i}.
\end{displaymath}
Assume we have a solution $z = (z_{1}, \ldots, z_{n})$ to the first $p$ consistency equations.  Then
\begin{displaymath}
\sum_{j=1}^{n}({r'_{ij} \llg(z_{j})} + {r''_{ij} \llg(1-z_{j})}) - c_{i} \pi {i} = 0 \qquad (i=1, \ldots, p).
\end{displaymath}
Thus,
\begin{displaymath}
\sum_{i=1}^{p}\lambda^{s}_{i}(\sum_{j=1}^{n}({r'_{ij} \llg(z_{j})} + {r''_{ij} \llg(1-z_{j})}) - c_{i} \pi {i}) = 0.
\end{displaymath}
Hence,
\begin{displaymath}
\sum_{j=1}^{n}\Big( (\sum_{i=1}^{p}\lambda^{s}_{i}r'_{ij}) \llg(z_{j}) +(\sum_{i=1}^{p}\lambda^{s}_{i}r''_{ij}) \llg(1-z_{j}) \Big) - (\sum_{i=1}^{p}\lambda^{s}_{i}c_{i}) \pi {i} = 0.
\end{displaymath}
This is the same as
\begin{displaymath}
\sum_{j=1}^{n}({r'_{sj} \llg(z_{j})} + {r''_{sj} \llg(1-z_{j})}) - c_{s} \pi {i} = 0.
\end{displaymath}
Therefore, the last $n - p$ consistency equations are determined by the first $p,$ so we only need the first $p$ equations to determine hyperbolicity.\hfill $\blacksquare$

In ~\cite{N-Z,B-P} it is proven that for a complete hyperbolic manifold, $\textrm{rank }\mathbf{R} = n - k$.  However, we need to prove hyperbolicity.  Neumann's work in Combinatorics of Triangulations and the Chern-Simons Invariant for Hyperbolic 3-Manifolds~\cite{N} tells us, without a priori knowledge of hyperbolicity, that $\textrm{rank }\mathbf{R} = n - k, \aand \mathbf{C}$ is determined by $\mathbf{R},$ so $\textrm{rank }\overline{\mathbf{R}} = n - k.$  This will be explained in Section~\ref{matrank}, ``Matrix Rank''.  Then, by the above proposition, we only need $n - k$ consistency equations to determine hyperbolicity.

\subsection{Cusp Conditions} \label{cc}
\begin{figure}
\begin{picture}(400,160)(-140,-35)
\put(0,0){\line(2,0){76}}
\put(76,0){\line(1,1){40}}
\put(116,40){\line(-3,2){78}}
\put(38,92){\line(-3,-2){78}}
\put(-40,40){\line(1,-1){40}}
\put(64,6){$w_{1}$}
\put(97,38){$w_{2}$}
\put(31,81){$w_{3}$}
\put(-36,38){$w_{4}$}
\put(-4,6){$w_{5}$}
\put(35,-6){$s_{1}$}
\put(99,17){$s_{2}$}
\put(85,63){$s_{3}$}
\put(-12,67){$s_{4}$}
\put(-32,19){$s_{5}$}
\put(180,40){\vector(-1,0){48}}
\put(182,40){$L_{\sigma_{22}}(v)$}
\put(180,-3){\vector(-2,1){60}}
\put(182,-3){$L_{\sigma_{21}}(v)$}
\put(85,90){\vector(-3,0){37}}
\put(130,86){\vector(-1,-2){16}}
\put(87,88){$L_{\sigma_{31}(v)}=L_{\sigma_{23}}(v)$}
\put(0,0){\line(-3,1){37}}
\put(0,0){\line(-1,-2){16}}
\put(76,0){\line(-1,-1){30}}
\put(76,0){\line(-1,-2){18}}
\put(76,0){\line(1,-2){16}}
\put(76,0){\line(3,-1){44}}
\put(76,0){\line(4,0){40}}
\put(116,40){\line(2,-1){28}}
\put(116,40){\line(1,2){15}}
\put(38,92){\line(4,1){40}}
\put(38,92){\line(1,2){15}}
\put(38,92){\line(-1,1){23}}
\put(-40,40){\line(-1,3){10}}
\put(-40,40){\line(-3,-2){25}}
\end{picture}
\caption{Simple Simplicial Loop, $\gamma,$ on Torus $T_{i}$} \label{loop}
\end{figure}
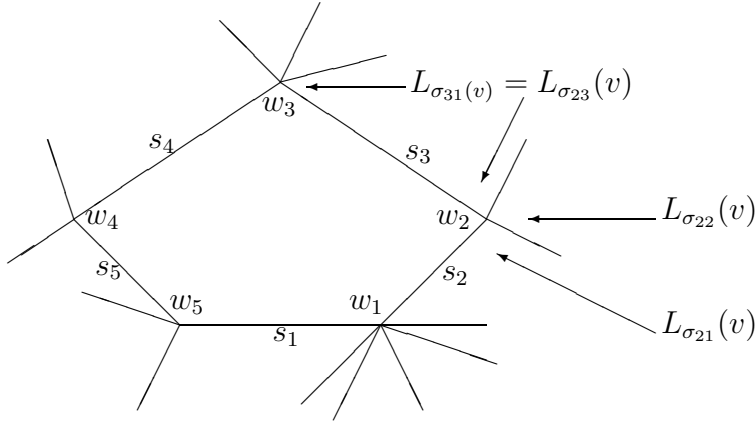

\begin{figure}
\begin{picture}(330,160)(-30,-35)
\put(76,0){\line(-6,1){68}}
\put(76,0){\line(1,1){40}}
\put(116,40){\line(3,1){57}}
\put(173,59){\line(3,-1){57}}
\put(230,40){\line(3,2){36}}
\put(266,64){\line(2,0){60}}
\put(64,6){$\tilde{w_{1}}$}
\put(102,42){$\tilde{w_{2}}$}
\put(170,63){$\tilde{w_{3}}$}
\put(220,45){$\tilde{w_{4}}$}
\put(264,68){$\tilde{w_{5}}$}
\put(35,9){$\tilde{s_{1}}$}
\put(85,21){$\tilde{s_{2}}$}
\put(140,54){$\tilde{s_{3}}$}
\put(197,53){$\tilde{s_{4}}$}
\put(241,57){$\tilde{s_{5}}$}
\put(296,70){$\tilde{s'_{1}}$}
\put(266,64){\line(1,-1){30}}
\put(266,64){\line(-1,-2){16}}
\put(76,0){\line(-1,-1){30}}
\put(76,0){\line(-1,-2){18}}
\put(76,0){\line(1,-2){16}}
\put(76,0){\line(3,-1){44}}
\put(76,0){\line(4,0){40}}
\put(116,40){\line(3,-1){32}}
\put(116,40){\line(0,-1){30}}
\put(173,59){\line(-3,-2){32}}
\put(173,59){\line(0,-1){30}}
\put(173,59){\line(3,-2){38}}
\put(230,40){\line(-1,-1){30}}
\put(230,40){\line(1,-1){30}}
\end{picture}
\caption{Developing Map Image of $\gamma$}  \label{loopdm}
\end{figure}
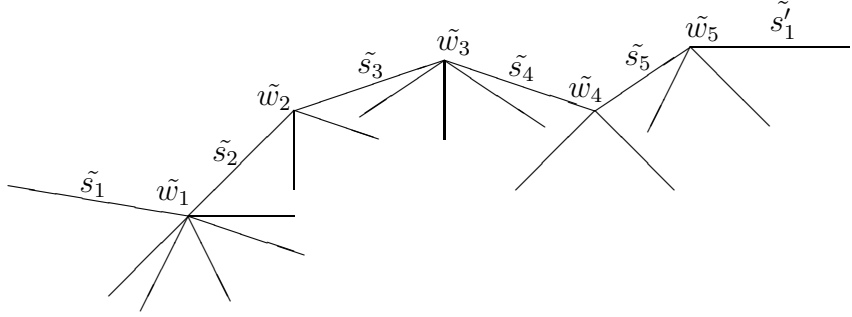

We now look at the $k$ cusps of $N.$  Details of the following discussion can be found in ~\cite{B-P}.  Let $T_{i}$ be the torus associated with the $i^{\mathrm{th}}$ cusp.  Select 2 simple oriented loops, $m_{i} \aand l_{i},$ on $T_{i},$ representing the 2 generators of the fundamental group of $T_{i}.$  Furthermore, $m_{i} \aand l_{i}$ can be chosen as simplicial loops with respect to $T_{i}\textrm{'s}$ triangulation.  Such a loop is composed of segments where each segment is an edge of some triangle $L_{\sigma_{q}}(v)\subset L(v) = T_{i}$, as identified earlier when describing the triangulation of $N$.  Let $\gamma$ be any simple simplicial oriented loop on $T_{i}$ consisting of $d$ segments, $s_{1}, \ldots,s_{d},$ and $d$ vertices, $w_{1}, \ldots, w_{d},$ where $w_{r}$ is the vertex at the end of $s_{r}$ as well as at the beginning of $s_{r+1}$ for $1 \leq r \leq d-1 \aand w_{d}$ is the vertex at the end of $s_{d}$ and beginning of $s_{1}$.  See Figure~\ref{loop}.
We lift  $\gamma$ to $\C = \R^{2},$ the universal cover of $T_{i},$ starting at the beginning of $s_{1}$ and map it to $\C$ by way of the developing map~\cite{Ratc,Thur}.  The resulting curve will consist of $d$ straight segments, $\tilde{s}_{1}, \ldots,\tilde{s}_{d},$ joined at the vertices $\tilde{w}_{r}$ for $1\leq r \leq d-1,$ as in $\gamma,$ except at $\tilde{w}_{d},$ which does not necessarily connect to the beginning of $\tilde{s}_{1}.$  So it starts at the beginning of $\tilde{s}_{1}$ and ends at the end of $\tilde{s}_{d}.$  Repeat the development map process, starting at the end of  $\tilde{s}_{d}$ and let $\tilde{s'}_{1}$ be the first segment this time, so  $\tilde{w}_{d}$ is the vertex between  $\tilde{s}_{d} \aand \tilde{s'}_{1}.$  See Figure~\ref{loopdm}.  Call this curve $\tilde{\gamma}.$
$\Aff(\C)$ can be regarded as $\C \rtimes \C^{*}$ with $(a,b)\in \C \rtimes \C^{*}$ such  that it represents $a + bx,$ an affine map of $\C$.  The \emph{dilation component of $(a,b)$} is $b$.  Thus, if an oriented triangle in $\C$ has two edges $\tilde{e}_{1} \aand \tilde{e}_{2}$ where $\tilde{e}_{1}$ ends in the vertex $\tilde{x},$ and $\tilde{e}_{2}$ begins at $\tilde{x},$ and the modulus of the triangle with respect to $\tilde{x}$ is $y,$ then the one and only orientation preserving similarity of $\C$ that takes $\tilde{e}_{1}$ to $\tilde{e}_{2}$ has dilation component equal to $-y$.  Remember, the modulus of the triangle with respect to $\tilde{x}$ is defined so that $\tilde{e}_{1}$ is identified with the edge from $0$ to $1$ and $\tilde{e}_{2}$ with the edge from $0$ to $y$ in the triangle with vertices $(0,1,y)$.  If $x_{r1}, \ldots, x_{rp_{r}}$ are the vertices of the $p_{r}$ triangles, $L_{\sigma_{r1}}(v),\ldots,L_{\sigma_{rp_{r}}}(v),$ that touch $\gamma$ at $w_{r},$ as in Figure~\ref{loop}, we get $p_{r}$ corresponding triangles, $\tilde{L}_{\sigma_{r1}}(v),\ldots,\tilde{L}_{\sigma_{rp_{r}}}(v),$ touching $\tilde{\gamma}$ at $\tilde{w}_{r}$ with $\tilde{x}_{r1},\ldots, \tilde{x}_{rp_{r}}$ the respective vertices of these triangles at $\tilde{w}_{r}$.  The ordering is such that $\tilde{s}_{r}$ is the first edge of $\tilde{L}_{\sigma_{r1}}(v),$ and $\tilde{s}_{r+1}$ is the second edge of $\tilde{L}_{\sigma_{rp_{r}}}(v),$ at $\tilde{w}_{r}$ unless $r=d,$ and then $\tilde{s'}_{1}$ is the second edge of $\tilde{L}_{\sigma_{dp_{d}}}(v).$ See Figure~\ref{loopdm}.   If the corresponding triangle modulii at $\tilde{w}_{r}$ are $y_{r1},\ldots,y_{rp_{r}},$ then the dilation component of the affine map that takes $\tilde{s}_{r}$ to $\tilde{s}_{r+1}$ is $- \prod_{i=1}^{p_{r}}{y_{ri}}$.  Orientation is responsible for the ``$-$" in the product.  Hence, the affine map that takes $\tilde{s}_{1}$ to $\tilde{s'}_{1}$ has dilation component of $\prod_{r=1}^{d} (-1)\prod_{i=1}^{p_{r}}{y_{ri}} = (-1)^{d}\prod_{r=1}^{d} \prod_{i=1}^{p_{r}}{y_{ri}}$.  Note that the modulus of $\tilde{L}_{\sigma_{ri}}(v)$ at $\tilde{x}_{ri}$ for $1\leq i \leq p_{r}$ is the same as the modulus of $L_{\sigma_{ri}}(v)$ at $x_{ri}$ for $1\leq i \leq p_{r},$  and this latter modulus has already been identified as either $z_{j},\, 1/(1 - z_{j}) \textrm{ or } 1 - 1/z_{j}$ for some $1\leq j \leq n$.  Therefore, the dilation component of the affine map that takes $\tilde{s}_{1}$ to $\tilde{s}'_{1}$ is of the form
\bdm
\pm 1\prod_{j=1}^{n} z_{j}^{\gamma'_{j}}(1-z_{j})^{\gamma''_{j}}.
\edm
The \emph{holonomy of the $(\C,\Aff(\C))$ structure} on $T_{i}$ is a map $\theta:\pi_{1}(T_{i})\to \Aff(\C)$ such that if [$\gamma$] is the element of $\pi_{1}(T_{i})$ represented by the loop $\gamma,$ then $\theta$ takes [$\gamma$] to the affine map that takes $\tilde{s_{1}}$ to $\tilde{s'_{1}}$.  This is a homomorphism that is well defined up to conjugacy class within $\Aff(\C)$.  However, any two elements of $\Aff(\C)$ within a conjugacy class have the same dilation component~\cite{B-P}, so the map

\begin{eqnarray*}
\psi_{i}: \pi_{1}(T_{i}) & \to & \C^* \qquad \textrm{such that}\\
{[}\gamma] & \to & \pm 1\prod_{j=1}^{n} z_{j}^{\gamma'_{ij}}(1-z_{j})^{\gamma''_{ij}}
\end{eqnarray*}
is a well defined homomorphism.  $\theta([\gamma])$ will be a translation if its dilation component is $1,$ so $\theta([\gamma])$ will be a translation when $\psi_{i}([\gamma]) = 1$. 

We now look at loops $m_{i} \aand l_{i}$.  For simplicity of notation, we also refer to the corresponding generators of $\pi_{1}(T_{i})$ as $m_{i} \aand l_{i}$ so 
\begin{eqnarray}
\psi_{i}(m_{i}) & = & \pm 1 \prod_{j=1}^{n} z_{j}^{m'_{ij}}(1-z_{j})^{m''_{ij}}\\
\psi_{i}(l_{i}) & = & \pm 1 \prod_{j=1}^{n} z_{j}^{l'_{ij}}(1-z_{j})^{l''_{ij}}
\end{eqnarray}
If the triangulation of $T_{i}$ causes $m_{i}$ to be a simplicial loop with $d$ segments and $d$ vertices, then its holonomy will be a non-trivial translation when $\psi_{i}(m_{i}) = 1$ and the sum of the arguments of the modulii at the $d$ vertices of $m_{i}$ is $d\pi$~\cite{B-P}.  Rewriting in log form, these requirements are expressed as
\bdm
\sum_{j=1}^{n}(m'_{ij} \llg(z_{j}) + m''_{ij} \llg(1-z_{j})) = c_{mi} \pi {i} \qquad\textrm{with }c_{mi}\in \Z.
\edm
Similarly, one can identify the log equation which sets the condition for the holonomy of $l_{i}$ to be a non-trivial translation.  It can be expressed as
\bdm
\sum_{j=1}^{n}(l'_{ij} \llg(z_{j}) + l''_{ij} \llg(1-z_{j})) = c_{li} \pi {i} \qquad\textrm{with }c_{li}\in \Z.
\edm
When the holonomy of the affine structure on $T_{i}$ has at least one non-trivial translation in its image, the affine structure is Euclidean~\cite{B-P}.  But a Euclidean structure on $T_{i}$ means that the $i^{\mathrm{th}}$ cusp is complete~\cite{B-P}, so the completeness equations for all of the $k$ cusps are
\beqn
\sum_{j=1}^{n}(m'_{ij} \llg(z_{j}) + m''_{ij} \llg(1-z_{j})) - c_{mi} \pi {i} = 0 \qquad (i=1, \ldots, k).
\eeqn

Now consider a hyperbolic manifold, $N,$ with $k$ cusps where $h$  of the cusps are complete, so the above completeness equations hold only for $k-h+1 \leq i \leq k$.  Let $T_{i}$ be a torus associated with one of the $k-h$ non-complete cusps.  If $p_{i} \aand q_{i}$ are co-prime integers, $(p_{i},q_{i})$ Dehn filling can be performed on this cusp.  In the literature, this process is frequently referred to as Dehn surgery, but it is really a filling.  In this case, $p_{i}m_{i} + q_{i}l_{i}$ is the generator of $\pi_{1}(T_{i})$ that is killed by Dehn filling.  In order to extend the hyperbolic structure on $N$ to the Dehn filling at this cusp, we need~\cite{N-Z,B-P}
\bdm
p_{i} \Big( \sum_{j=1}^{n}(m'_{ij} \llg(z_{j}) + m''_{ij} \llg(1-z_{j})) - c_{mi} \pi {i} \Big) + 
\edm
\beqn
\qquad \qquad q_{i} \Big( \sum_{j=1}^{n}(l'_{ij} \llg(z_{j}) + l''_{ij} \llg(1-z_{j})) - c_{li} \pi {i} \Big) = 2\pi {i}.
\eeqn
That is
\bdm
\sum_{j=1}^{n}\Big( (p_{i}m'_{ij} + q_{i}l'_{ij}) \llg(z_{j}) + (p_{i}m''_{ij} + q_{i}l''_{ij}) \llg(1-z_{j}) \Big) = c_{si}\pi {i} \qquad\textrm{with }c_{si}\in \Z.
\edm
Therefore, if the equations
\beqn
\sum_{j=1}^{n}\Big( (p_{i}m'_{ij} + q_{i}l'_{ij}) \llg(z_{j}) + (p_{i}m''_{ij} + q_{i}l''_{ij}) \llg(1-z_{j}) \Big) = c_{si}\pi {i}
\eeqn
\bdm
 \qquad \qquad \qquad \qquad \qquad \qquad \qquad \qquad \qquad \qquad \qquad \qquad (i =  1, \ldots, k-h)
\edm
are satisfied, $M,$ the manifold derived from $N$ by Dehn filling on the $k-h$ cusps, will be hyperbolic near these cusps.

The last step in identifying the equations is the selection of the appropriate $n - k$ consistency equations.  Let $s'_{ij} = p_{i}m'_{ij} + q_{i}l'_{ij} \aand s''_{ij} = p_{i}m''_{ij} + q_{i}l''_{ij},$ and define the matrices $\mathbf{M}, \, \mathbf{L}, \, \mathbf{S} \aand \mathbf{M_{h}}$ as $\mathbf{R}$ is defined on page~\pageref{Rmat} so that
\begin{displaymath}
\mathbf{M} =
\left( \begin{array}{cccccc}
m'_{11} & \ldots & m'_{1n} & m''_{11} & \ldots & m''_{1n} \\
\vdots & \ddots & \vdots & \vdots & \ddots & \vdots \\
m'_{(k)1} & \ldots & m'_{(k)n} & m''_{(k)1} & \ldots & m''_{(k)n} 
\end{array} \right)  \label{Mmat}
\edm

\bdm
\mathbf{L} =
\left( \begin{array}{cccccc}
l'_{11} & \ldots & l'_{1n} & l''_{11} & \ldots & l''_{1n} \\
\vdots & \ddots & \vdots & \vdots & \ddots & \vdots \\
l'_{(k)1} & \ldots & l'_{(k)n} & l''_{(k)1} & \ldots & l''_{(k)n} 
\end{array} \right)  \label{Lmat}
\edm

\begin{displaymath}
\mathbf{S} =
\left( \begin{array}{cccccc}
s'_{11} & \ldots & s'_{1n} & s''_{11} & \ldots & s''_{1n} \\
\vdots & \ddots & \vdots & \vdots & \ddots & \vdots \\
s'_{(k-h)1} & \ldots & s'_{(k-h)n} & s''_{(k-h)1} & \ldots & s''_{(k-h)n} 
\end{array} \right)   \label{Smat} 
\end{displaymath}

\begin{displaymath}
\mathbf{M_{h}} =
\left( \begin{array}{cccccc}
m'_{(k-h+1)1} & \ldots & m'_{(k-h+1)n} & m''_{(k-h+1)1} & \ldots & m''_{(k-h+1)n} \\
\vdots & \ddots & \vdots & \vdots & \ddots & \vdots \\
m'_{(k)1} & \ldots & m'_{(k)n} & m''_{(k)1} & \ldots & m''_{(k)n} 
\end{array} \right)
\edm
Let
\bdm
\mathbf{U} =
\left( \begin{array}{c}
\mathbf{\mathbf{S}}\\
\mathbf{\mathbf{M_{h}}}
\end{array} \right )   \label{Umat}
\edm
We will see that rank $\mathbf{U} = k$.  We can select $n - k$ consistency equations so that their rows in $\mathbf{R}$ are linearly independent, and when concatenated with $\mathbf{U},$ give an $n \times (2n)$ matrix of rank $n$.  The reasons for this are a consequence of~\cite{N}, and will be explained in Section~\ref{matrank}, ``Matrix Rank''.  We will assume, without loss of generality, that the last $n-k$ out of $n$ consistency equations are the ones we want.

In summary, we have\\
$n-k$ consistency equations,
\bdm
\sum_{j=1}^{n}(r'_{ij} \llg(z_{j}) + r''_{ij} \llg(1-z_{j})) - c_{i} \pi {i} = 0 \qquad  (i=k+1, \ldots, n),
\edm
$k-h$ surgery equations,
\bdm
\sum_{j=1}^{n}\Big( (p_{i}m'_{ij} + q_{i}l'_{ij}) \llg(z_{j}) + (p_{i}m''_{ij} + q_{i}l''_{ij}) \llg(1-z_{j}) \Big) - c_{si} \pi {i} =  0
\edm
\bdm
 \qquad \qquad \qquad \qquad \qquad \qquad \qquad \qquad \qquad \qquad \qquad \qquad (i =  1, \ldots, k-h),
\edm
and $h$ completeness equations,
\bdm
\sum_{j=1}^{n}(m'_{ij} \llg(z_{j}) + m''_{ij} \llg(1-z_{j})) - c_{mi} \pi {i} = 0 \qquad (i=k-h+1, \ldots, k).
\edm
giving a total of $n$ equations that must have a simultaneous solution to make a manifold complete hyperbolic.

\subsection{Matrix Rank} \label{matrank}
In~\cite{N}, Neumann has constructed a chain complex, $\J,$ and described its homology.  Using the terminology of Section~\ref{intro}, ``Introduction,'' with respect to the triangulation of $N \aand M,$ let $K$ be the gluing of the $n$ tetrahedra, $\sigma_{1}, \ldots, \sigma_{n}$.  The modules of the chain complex are $C_{0},\, C_{1} \aand J,$ where
\begin{enumerate}
\item $C_{0} = \Z$ module generated by the $k$ vertices of $K$.  Each vertex will be associated with a cusp of $N,$ and the torus that is the link of the vertex.
\item $C_{1} = \Z$ module generated by $E_{1}, \ldots, E_{n},$ the $n$ edges of $K$.
\item With regard to $J,$ for each tetrahedron, $\sigma_{j},$ label the edges as $e_{j1}, \ldots, e_{j6}$ according to the associated parameters as:
\bdm
\begin{array}{lll}
e_{j1}=z & \qquad e_{j2}=\frac{1}{1-z} & \qquad e_{j3}=1-\frac{1}{z}\\
e_{j4}=z & \qquad e_{j5}=\frac{1}{1-z} & \qquad e_{j6}=1-\frac{1}{z}
\end{array}
\edm
Let $\Jsig{j}=\Z$ module generated by the six edges of $\sigma_{j}$ with the relations $e_{j\tau} - e_{j(\tau +3)} = 0 \textrm{ for } 1 \leq \tau \leq 3 \aand e_{j1}+e_{j2}+e_{j3} = 0$.  Thus, opposite edges of the tetrahedron are represented by the same element of $\Jsig{j}, \aand e_{j3}$ can be defined in terms of $e_{j1} \aand e_{j2}$.  This means that $e_{j1} \aand e_{j2}$ generate the $\Z$ module, $\Jsig{j}$.  Let
\bdm
J = \coprod_{1 \leq j \leq n} \Jsig{j}.
\edm
\end{enumerate}
The chain complex sequence is
\bdm
\J: \quad 0 \rightarrow C_{0} \stackrel{\alpha}{\rightarrow} C_{1} \stackrel{\beta}{\rightarrow} J \stackrel{\beta^{*}}{\rightarrow} C_{1} \stackrel{\alpha^{*}}{\rightarrow} C_{0} \rightarrow 0.
\edm
We have $\alpha, \, \beta, \, \alpha^{*} \aand \beta^{*}$ defined as follows:
\begin{enumerate}
\item $\alpha:C_{0} \to C_{1},$ where $\alpha$ takes a vertex to the sum of the edges containing the vertex, with an edge counted twice if both ends of the edge are at the vertex.
\item $\beta:C_{1} \to J$ can be defined by letting 
\bdm
E_{i} \to \sum_{1 \leq j \leq n} \sum_{\stackrel{1 \leq \tau \leq 6}{\scriptscriptstyle E_{i}\, \mathrm{is}\, \mathrm{identified}\, \mathrm{with}\, e_{j\tau}}}e_{j\tau}
\edm
We have the sum $\sum_{\stackrel{1 \leq \tau \leq 6}{\scriptscriptstyle E_{i}\, \mathrm{is}\, \mathrm{identified}\, \mathrm{with}\, e_{j\tau}}}e_{j\tau} \in \Jsig{j}$ because more than one edge of $\sigma_{j}$ can be identified with $E_{i}$.
\item To define $\beta^{*}:J \to C_{1},$ note that for each $\sigma_{j},$ we have the edge set $\{ e_{j1}, \ldots, e_{j6} \}$.  Let $\rho:\{ e_{j1}, \ldots, e_{j6} \} \to \{ E_{1}, \ldots ,E_{n} \}$ be such that $\rho(e_{j\tau}) = E_{i}$ when $e_{j\tau}$ is identified with the edge $E_{i}$.  Then, let
\bdm
\beta^{*}(e_{j\tau}) = \rho(e_{j(\tau + 1)}) - \rho(e_{j(\tau + 2)}) + \rho(e_{j(\tau + 4)}) - \rho(e_{j(\tau + 5)}) \quad \textrm{(indices mod 6)}
\edm  
That is, $\beta^{*}$ takes $e_{j\tau}$ to the alternating sum of the edges of $N$ identified with the edges of $\sigma_{j}$ that touch $e_{j\tau}$.
\item $\alpha^{*}:C_{1} \to C_{0},$ where $\alpha^{*}$ sends an edge, $E_{i},$ to the sum of its end points.
\end{enumerate}
$N$ is the interior of a compact manifold, $\overline{N},$ whose boundary is the union of the $k$ tori, $T_{1}, \ldots , T_{k},$ that are the links of the vertices of $K$. 
\begin{lemma} \label{seq}
When tensored with $\Q,$ the sequence, $\J,$ is exact except in the middle, where its homology is $H_{1}(\partial \overline{N};\Q)= \coprod_{1 \leq i \leq k} H_{1}(T_{i};\Q)$.
\end{lemma}
For a proof, see~\cite{N}.  We use this to compute the rank of $\mathbf{R}$.  However, we will use the original chain with coefficients in $\Z$ to show that the rank of the matrix obtained by concatenating $\mathbf{U},$ as defined on page~\pageref{Umat}, with $n-k$ linearly independent rows of $\mathbf{R},$ is $n$.
\subsubsection{Rank of $\Rm$} 
The matrix of the linear transformation, $\beta,$ is closely related to $\Rm^{t},$ the transpose of $\Rm,$ and they have the same rank.  Since $\textrm{rank}\, \Rm = \textrm{rank}\, \Rm^{t}, \ \textrm{rank}\, \Rm = \textrm{rank}$ of the matrix of $\beta$.  The edges $E_{1}, \ldots, E_{n}$ are a basis of $C_{1}$ as a vector space, so the vectors $\beta(E_{i}) \textrm{ for } 1 \leq i \leq n$ are the columns of the matrix of $\beta$.  From the definition of $\beta,$ we see that in $\Jsig{j},$ 
\beqn
\beta(E_{i}) = \sum_{\stackrel{1 \leq \tau \leq 6}{\scriptscriptstyle E_{i}\, \mathrm{is}\, \mathrm{identified}\, \mathrm{with}\, e_{j\tau}}}e_{j\tau}  \qquad \textrm{ modulo relations on }J.
\eeqn                      
Thus, if:\\
\hspace*{.2in}$e_{j1}$ or $e_{j4}$ occur, it means $E_{i}$ is identified with the $z_{j}$ parameter \\
\hspace*{.2in}$e_{j2}$ or $e_{j5}$ occur, it means $E_{i}$ is identified with the $\frac{1}{1-z_{j}}$ parameter \\
\hspace*{.2in}$e_{j3}$ or $e_{j6}$ occur, it means $E_{i}$ is identified with the $1 - \frac{1}{z_{j}}$ parameter.\\
In $\Jsig{j}, \, e_{j3} = -e_{j1} - e_{j2}$; also, $1 - \frac{1}{z_{j}} =  -\frac{1-z_{j}}{z_{j}}$.  Hence, the sum of the coefficients of $e_{j1}$ in $\beta(E_{i})$ is $r'_{ij},$ the sum of the exponents of $z_{j}$ with respect to the edge $E_{i}$ in the consistency equations, and the sum of the coefficients of $e_{j2}$ in $\beta(E_{i})$ is $-r''_{ij},$ which is $-1$ times the sum of the exponents of $1-z_{j}$ with respect to the edge $E_{i}$ in the consistency equations, as seen on page~\pageref{rval}.  Consequently, $\Rmu,$ the $2n \times n$ matrix of $\beta$ is
\begin{displaymath}
\Rmu =
\left( \begin{array}{cccccc}
r'_{11} & \ldots & r'_{n1}\\
-r''_{11} & \ldots & -r''_{n1} \\
\vdots & \ddots & \vdots \\
r'_{1n} & \ldots & r'_{nn}\\
-r''_{1n} & \ldots & -r''_{nn} 
\end{array} \right) 
\end{displaymath}
We see that $\textrm{rank}\, \Rmu = \textrm{rank}\, \Rm^{t},$ so $\textrm{rank}\, \Rm = \textrm{rank}\, \Rmu$.  By definition, the rank of $\Rmu$ is equal to the dimension of the image of $\beta$.  By Lemma~\ref{seq}, $\alpha$ is injective, making $\textrm{dim}\,\textrm{im}(\alpha) = \textrm{dim}\,C_{0} = k, \aand \im(\alpha) = \ker(\beta)$.  The matrix of $\beta$ this way would still be $\Rmu,$ so
\begin{eqnarray*}
\textrm{rank}\, \Rmu & = & \textrm{dim}\,\textrm{im}(\beta)\\
& = & \textrm{dim}\, C_{1} - \textrm{dim}\, \textrm{kernel}(\beta)\\
& = & \textrm{dim}\, C_{1} - \textrm{dim}\, \textrm{im}(\alpha)\\
& = & n - k
\end{eqnarray*}    
Therefore, $\textrm{rank}\, \Rm = n - k$.  Let $\Rm, \, \Cm \aand \Rmo$ be the matrices associated with the consistency equations, as on page~\pageref{Rmat}.  Consider the matrix equation $\Rm \cdot x = \mathbf{-}\Cm$.  In~\cite{N} it is proved that there is an $\tilde{x} \in \Q^{2n}$ that is a solution.  Then $\mathbf{-}\Cm$ is a linear combination of the columns of $\Rm,$ so $\Rm$ concatenated with $\mathbf{-}\Cm$ has the same rank as $\Rm$ since row rank is the same as column rank.  That is, 
\begin{eqnarray*}
n - k & = & \textrm{rank}\, \Rm = \textrm{column rank}\, \Rm = \textrm{column rank}\, (\Rm|\mathbf{-}\Cm)\\ 
& = &\textrm{column rank}\, (\Rm|\Cm)=\textrm{column rank}\, \Rmo = \textrm{rank}\, \Rmo,
\end{eqnarray*}
so $\textrm{rank}\, \Rmo = n - k$.  Let $\Rb = $ matrix consisting of $n - k$ linearly independent rows of $\Rm$.
\subsubsection{Rank of $(\Sm|\Mhm|\Rb)$}
For now we will include all $k$ cusps of $N$.  Let $S_{1}(\partial \overline{N}) = \Z$ module of simplicial 1-chains, $Z_{1}(\partial \overline{N}) = \Z$ module of 1-cycles and $B_{1}(\partial \overline{N}) = \Z$ module of 1-boundaries.  Let $e_{j\tau} \in \Jsig{j}$ for $\tau = 1,2$.  If the two vertices at the ends of $e_{j\tau}$ in $\sigma_{j}$ are $v_{j\tau 1} \aand v_{j\tau 2},$ let $\zeta_{j\tau 1} \aand \zeta_{j\tau 2}$ be the respective edges of $\Lsig{j}(v_{j\tau 1}) \aand \Lsig{j}(v_{j\tau 2})$ that do not intersect $e_{j\tau}$.  Do the same for $e_{j(\tau+3)},$ so we have four 1-simplices identified in $\partial \overline{N}$.  They are $\zeta_{j\tau 1},$ $\zeta_{j\tau 2},$ $\zeta_{j(\tau +3) 1}$ and $\zeta_{j(\tau + 3) 2},$ with one for each vertex of $\sigma_{j}$.  Now define $\hgamma_{0}$.
\bqr
\hgamma_{0}: \Jsig{j} & \to & S_{1}(\partial \overline{N})\\
e_{j\tau} & \to & \zeta_{j\tau 1} + \zeta_{j\tau 2} + \zeta_{j(\tau +3) 1} + \zeta_{j(\tau + 3) 2}
\eqr
We have, by~\cite{N},  
\bqr
\hgamma_{0}: \textrm{ im}(\beta) & \to & B_{1}(\partial \overline{N})\\
\hgamma_{0}: \textrm{ ker}(\beta^{*}) & \to & Z_{1}(\partial \overline{N})
\eqr
so there is the induced map
\bdm
\hgamma: \textrm{ker}(\beta^{*})/\textrm{im}(\beta) \to H_{1}(\partial \overline{N}) = \coprod_{1 \leq i \leq k}H_{1}(T_{i})
\edm
Next, let $\hdelta_{0}:\, H_{1}(\partial \overline{N}) \to J$ be defined as follows.  Let $\Gamma$ be a simple simplicial loop on the torus, $T_{i},$ associated with the $i^{th}$ cusp of $N$.  In figure~\ref{loop}, $\gamma$ is such a loop.  Each vertex, $w_{r},$ of $\gamma,$ is the vertex of $p_{r}$ triangles $\Lsigv{r_{1}}, \ldots , \Lsigv{r_{p_{r}}}$ where $T_{i}$ is the link of $v,$ a vertex of $K$.  Define the simple cellular path $\overline{\Gamma},$ by starting at the midpoint of the edge of $\Lsigv{11}$ that ends in $w_{1}$ but is not $s_{1}$.  Continue across the $\{\Lsigv{1q}\}_{2 \leq q \leq p_{1}-1}$ by crossing from one triangle to another at the midpoint of the edges that have $w_{1}$ as a vertex, ending at the edge of $\Lsigv{1p_{1}}$ that is not $s_{2}$.  Then continue across $\Lsigv{1p_{1}} = \Lsigv{21}$ to the edge of $\Lsigv{21}$ that has $w_{2}$ as a vertex but is not $s_{2}$.  Repeat the process until the loop is closed by going from the edge of $\Lsigv{dp_{d}} = \Lsigv{11}$ that contains $w_{d}$ but is not $s_{1}$ to the starting point.  When $\overline{\Gamma}$ crosses $\Lsigv{rq}$ for $2 \leq q \leq p_{r}-1,$ it goes counterclockwise around the vertex $w_{r},$ as a vertex of $\Lsigv{rq},$ and when it crosses $\Lsigv{rp_{r}} = \Lsigv{(r+1)1},$ it goes clockwise around the vertex of this triangle that is opposite to $s_{r+1}$.  When one of these vertices belongs to the triangle $\Lsigv{rq},$ the vertex is associated with an edge, $e_{rq\tau},$ of $\sigma_{rq}$ for some $1 \leq \tau \leq 6,$ as defined at the beginning of Section~\ref{matrank}, ``Matrix Rank'', and this edge is an element of $\Jsig{rq} \subset J$.  To each of these edges assign a ``$+$'' if $\overline{\Gamma}$ goes around its corresponding vertex counterclockwise, and a ``$-$'' if $\overline{\Gamma}$ goes around its corresponding vertex clockwise.  $\Gamma$ is homotopic to $\overline{\Gamma},$ so we can define $\hdelta_{0}:\, Z_{1}(\partial \overline{N}) \to J$ such that $\hdelta_{0}(\overline{\Gamma}) = \hdelta_{0}(\Gamma)$ is the signed sum of these edges in $J$.  That is,
\beqn
\hdelta_{0}(\overline{\Gamma}) = \sum_{\stackrel{1 \leq r \leq d}{\scriptscriptstyle 2 \leq q \leq p_{r}}}(-1)^{t}e_{rq\tau}
\eeqn  
where $t = 0$ when $e_{rq\tau}$ is assigned a ``$+$'' and $t=1$ when $e_{rq\tau}$ is assigned a ``$-$''.
In $\Jsig{rq}, \, e_{rq\tau} = e_{rq(\tau + 3)}\textrm{ for }1\leq \tau \leq 3$ with the last subscript mod $6,$ and $e_{rq1} + e_{rq2} + e_{rq3} = 0,$ so $-e_{rq\tau} = e_{rq(\tau + 1)} + e_{rq(\tau + 2)}$ with the last two subscripts mod $6$.  Therefore, when $e_{rq\tau}$ is assigned a ``$-$'', we substitute $e_{rq(\tau + 1)} + e_{rq(\tau + 2)}$ with both subscripts mod $6$. Hence,
\beqn
\hdelta_{0}(\overline{\Gamma}) = \sum_{\stackrel{1 \leq r \leq d}{\scriptscriptstyle 1 \leq q \leq p_{r}}}e_{rq\tau} \label{gam}
\eeqn
where $e_{rq\tau}$ is an edge of $\sigma_{rq}$ that is associated with $w_{r},$ a vertex of $\Lsigv{rq}$ and $w_{r}$ is a vertex of the simple simplicial loop $\Gamma$ in $T_{i}$.  The relations of $J$ also mean that $e_{rq3} = e_{rq6} = -e_{rq1} - e_{rq2},$ so
\beqn
\hdelta_{0}(\overline{\Gamma}) = \sum_{\stackrel{1 \leq j \leq n}{\scriptscriptstyle \overline{\Gamma} \, \mathrm{crosses}\,\Lsigv{j}}}g'_{j\overline{\Gamma}}e_{j1} + g''_{j\overline{\Gamma}}e_{j2}
\eeqn
where, with respect to $\sigma_{j}, \, g'_{j\overline{\Gamma}}$ is the number of occurrences of the $z_{j}$ parameter minus the number of occurrences of the $1 - \frac{1}{z_{j}}$ parameter and $g''_{j\overline{\Gamma}}$ is the number of occurrences of the $\frac{1}{1 - z_{j}}$ parameter minus the number of occurrences of the $1 - \frac{1}{z_{j}}$ parameter in Equation~\ref{gam}.  

Now let $m_{i} \aand l_{i} \textrm{ for } 1 \leq i \leq k$ be the meridianal and longitudinal simple simplicial loops on $T_{i},$ as in Section~\ref{cc}, ``Cusp Conditions''.  We get corresponding $\overline{m}_{i} \aand \overline{l}_{i},$ constructed as $\overline{\Gamma}$ was, where $m_{i} \aand l_{i}$ are homologous to $\overline{m}_{i} \aand \overline{l}_{i},$ respectively.  So $\overline{m}_{i} \aand \overline{l}_{i}$ are the generators of $H_{1}(T_{i})$ and their image under $\hdelta_{0}$ are two columns of $\MLu,$ the matrix of $\hdelta_{0}$.  These two columns are of the form
\bqr
\vec{g}_{m_{i}} & = & (g'_{1\overline{m}_{i}},g''_{1\overline{m}_{i}}, \ldots, g'_{n\overline{m}_{i}},g''_{n\overline{m}_{i}})\qquad \textrm{with } g'_{j\overline{m}_{i}} = m'_{ij} \aand g''_{j\overline{m}_{i}} = -m''_{ij}\\      
\vec{g}_{l_{i}} & = & (g'_{1\overline{l}_{i}},g''_{1\overline{l}_{i}}, \ldots, g'_{n\overline{}_{i}},g''_{n\overline{l}_{i}})\qquad \textrm{with } g'_{j\overline{l}_{i}} = l'_{ij} \aand g''_{j\overline{l}_{i}} = -l''_{ij}
\eqr 
where $m'_{ij},m''_{ij} \aand l'_{ij},l''_{ij}$ are the components of the matrices $\Mm \aand \Lm$ from ``Cusp Conditions'' on page~\pageref{Mmat}.  Let $\MLm$ be $\Mm$ concatenated with $\Lm$.  For each generator of $H_{1}(\partial \overline{N}) = \coprod_{1 \leq i \leq k}H_{1}(T_{i}),$ there is a column in the matrix of $\hdelta_{0},$ so $\MLu$ has $2k$ columns and $2n$ rows, where the $(2j - 1)^{th}$ row of $\MLu$ is equal to the $j^{th}$ column of $\MLm$ and the $2j^{th}$ row of $\MLu$ is $(-1)$ times the $(n+j)^{th}$ column of $\MLm$.  Thus, $\textrm{rank}\, \MLu = \textrm{rank}\, \MLu^{t} = \textrm{rank}\, \MLm$.  The next step is to show that $\textrm{rank}\, \MLu = 2k$.  We have $\im(\hdelta_{0}) \subset \ker(\beta^{*}),$ with $\hdelta_{0}(B_{1}(\partial\overline{N})) \subset \im(\beta),$ so there is the induced map
\bdm
\hdelta:\,H_{1}(\partial \overline{N}) \to \ker(\beta^{*})/\im(\beta)
\edm
Now $\hgamma \hdelta: \,H_{1}(\partial \overline{N}) \to H_{1}(\partial \overline{N})$ is multiplication by $2$~\cite{N}, so $\hdelta_{0}$ must be injective.  Consequently, the matrix of $\hdelta_{0}$ has maximal rank, which is $2k,$ making the $2k$ vectors, $\{ \vec{g}_{m_{i}}, \, \vec{g}_{l_{i}} \}_{1 \leq i \leq k},$ linearly independent.

$M$ is derived from $N$ by the Dehn filling of $k - h$ cusps of $N$ with filling coefficients of $(p_{i},q_{i})$ for $1 \leq i \leq k - h$.  Let $\vec{g}_{s_{i}} = p_{i}\vec{g}_{m_{i}} + q_{i}\vec{g}_{l_{i}}$ for $1 \leq i \leq k - h$.
\begin{lemma} \label{l2}
The $k + h$ vectors
\bdm
\{ \vec{g}_{s_{1}}, \ldots, \vec{g}_{s_{k-h}}, \vec{g}_{m_{k-h+1}}, \ldots, \vec{g}_{m_{k}}, \vec{g}_{l_{k-h+1}}, \ldots, \vec{g}_{l_{k}} \}
\edm
are linearly independent.
\end{lemma}

\emph{Proof.} 
Assume otherwise.  Then there exists $\phi_{si}$ for $1 \leq i \leq k - h$ and $\xi_{mi} \aand \varphi_{li}$ for $k - h +1 \leq i \leq k$ such that
\begin{eqnarray*}
0 & = & \sum_{1 \leq i \leq k-h}\phi_{si} \vec{g}_{s_{i}} + \sum_{k-h+1 \leq i \leq k}(\xi_{mi}\vec{g}_{m_{i}} + \varphi_{li}\vec{g}_{l_{i}})\\
& = &  \sum_{1 \leq i \leq k-h}\phi_{si}(p_{i}\vec{g}_{m_{i}} + q_{i}\vec{g}_{l_{i}} ) + \sum_{k-h+1 \leq i \leq k}(\xi_{mi}\vec{g}_{m_{i}} + \varphi_{li}\vec{g}_{l_{i}}) \\            
& = &  \sum_{1 \leq i \leq k-h}\phi_{si}p_{i}\vec{g}_{m_{i}} + \sum_{1 \leq i \leq k-h}\phi_{si}q_{i}\vec{g}_{l_{i}} + \sum_{k-h+1 \leq i \leq k}(\xi_{mi}\vec{g}_{m_{i}} + \varphi_{li}\vec{g}_{l_{i}})            
\end{eqnarray*}
We have just seen that $\{ \vec{g}_{m_{i}}, \, \vec{g}_{l_{i}} \}_{1 \leq i \leq k},$ is linearly independent, so $\xi_{mi} = \varphi_{li} = 0$ for $k-h+1 \leq i \leq k$ and $\phi_{si}p_{i} = \phi_{si}q_{i} = 0$ for $1 \leq i \leq k - h$.  But at least one of $p_{i}$ or $q_{i}$ is not $0,$ so $\phi_{si} = 0$ for $1 \leq i \leq k-h$. \hfill $\blacksquare$

Since $\textrm{rank}\,\Rmu = n - k,$ select $n - k$ linearly independent vectors in $\im(\beta)$ that are columns of the matrix $\Rmu,$ and denote them by $\vec{g}_{\beta_{i}}$ for $k + 1 \leq i \leq n$.  Observe that $\im(\hdelta_{0}) \cap \im(\beta) = \{ 0 \},$ because otherwise, there is a non-trivial $x \in H_{1}(\partial \overline{N})$ such that $\hgamma_{0} \hdelta_{0}(x) = \hgamma_{0}(\textrm{element of im}(\beta)) \in B_{1}(\partial \overline{N})$.  Then $\hgamma \hdelta (x) = 0 \in H_{1}(\partial \overline{N})$.  But $\hgamma \hdelta$ is multiplication by $2$ on $H_{1}(\partial \overline{N}),$ so $x = 0,$ which is a contradiction.
\begin{lemma} \label{l3}
Let\\
1) $\Smu =$ the $2n \times (k - h)$ matrix whose columns are the vectors $\vec{g}_{s_{i}},$ for $1 \leq i \leq k-h$\\
2) $\Mhmu =$ the $2n \times h$ matrix whose columns are the linearly independent vectors $\vec{g}_{m_{i}},$ for $k-h + 1 \leq i \leq k$\\
3) $\Rmb =$ the $2n \times (n-k)$ matrix whose columns are the linearly independent vectors $\vec{g}_{\beta_{i}},$ for $k + 1 \leq i \leq n$\\
Concatenate these matrices to get the $2n \times n$ matrix $\Fmu = (\Smu|\Mhmu|\Rmb)$.  Then $\textrm{Rank} \, \Fmu = n$. 
\end{lemma}

\emph{Proof.}
Assume otherwise.  Then the vectors that are the columns of $\Fmu$ are not linearly independent, so there are $\xi_{si}$ for $1 \leq i \leq k - h$, $\varphi_{mi}$ for $k - h +1 \leq i \leq k$ and $\phi_{\beta i}$ for $k + 1 \leq i \leq n,$ where not all are zero, such that
\beqn
0 = \sum_{1 \leq i \leq k-h}\xi_{si}\vec{g}_{s_{i}} + \sum_{k-h+1 \leq i \leq k}\varphi_{mi}\vec{g}_{m_{i}} + \sum_{k+1 \leq i \leq n}\phi_{\beta i} \vec{g}_{\beta_{i}} \label{all}
\eeqn
Therefore,
\begin{eqnarray*}
0 & = & \hgamma_{0}(0)\\ 
& = & \hgamma_{0}\Big( \sum_{1 \leq i \leq k-h}\xi_{si}\vec{g}_{s_{i}} + \sum_{k-h+1 \leq i \leq k}\varphi_{mi}\vec{g}_{m_{i}} + \sum_{k+1 \leq i \leq n}\phi_{\beta i} \vec{g}_{\beta_{i}} \Big)\\
& = &  \hgamma_{0}\Big( \sum_{1 \leq i \leq k-h}\xi_{si}\vec{g}_{s_{i}} + \sum_{k-h+1 \leq i \leq k}\varphi_{mi}\vec{g}_{m_{i}} \Big) + \hgamma_{0}\Big( \sum_{k+1 \leq i \leq n}\phi_{\beta i} \vec{g}_{\beta_{i}} \Big)\\
& = &  \hgamma_{0}\Big[\hdelta_{0}\Big( \sum_{1 \leq i \leq k-h}\xi_{si}(p_{i}\overline{m}_{i} + q_{i}\overline{l}_{i}) + \sum_{k-h+1 \leq i \leq k}\varphi_{mi}\overline{m}_{i} \Big)\Big] + \hgamma_{0}(\textrm{element in im}(\beta))
\end{eqnarray*}
But $\hgamma_{0}(\im(\beta)) \subset B_{1}(\partial \overline{N}),$ so $\hgamma \hdelta \Big( \sum_{1 \leq i \leq k-h}\xi_{si}(p_{i}\overline{m}_{i} + q_{i}\overline{l}_{i}) + \sum_{k-h+1 \leq i \leq k}\varphi_{mi}\overline{m}_{i} \Big) = 0$.  Therefore, $\sum_{1 \leq i \leq k-h}\xi_{si}(p_{i}\overline{m}_{i} + q_{i}\overline{l}_{i}) + \sum_{k-h+1 \leq i \leq k}\varphi_{mi}\overline{m}_{i} = 0$ since $\hgamma \hdelta$ is injective.  Hence, 
\begin{eqnarray*}
0 & = & \hdelta_{0}(0)\\
& = & \hdelta_{0}\Big( \sum_{1 \leq i \leq k-h}\xi_{si}(p_{i}\overline{m}_{i} + q_{i}\overline{l}_{i}) + \sum_{k-h+1 \leq i \leq k}\varphi_{mi}\overline{m}_{i} \Big)\\
& = & \sum_{1 \leq i \leq k-h}\xi_{si}\hdelta_{0}(p_{i}\overline{m}_{i} + q_{i}\overline{l}_{i}) + \sum_{k-h+1 \leq i \leq k}\varphi_{mi}\hdelta_{0}(\overline{m}_{i})\\ 
& = & \sum_{1 \leq i \leq k-h}\xi_{si}\vec{g}_{s_{i}} + \sum_{k-h+1 \leq i \leq k}\varphi_{mi}\vec{g}_{m_{i}}
\end{eqnarray*}
By Lemma~\ref{l2}, $\xi_{si}$ for $1 \leq i \leq k - h \aand \varphi_{mi}$ for $k - h +1 \leq i \leq k$ are all zero.  Then, Equation~\ref{all} becomes $0 = \sum_{k + 1 \leq i \leq n}\phi_{\beta i}\vec{g}_{\beta_{i}}$.  However, the $\vec{g}_{\beta_{i}}$ for $k+1 \leq i \leq n$ were selected to be linearly independent, so $\phi_{\beta i} = 0$ for $k+1 \leq i \leq n$.  This is a contradiction. \hfill $\blacksquare$ 
\begin{corollary}
Each column of $\Rmb$ has a corresponding row in $\Rm,$ the matrix associated with the consistency equations.  Let $\Rb$ be the matrix comprised of only these $n-k$ rows of $\Rm$ and let 
\bdm
\Fm = \left( \begin{array}{c}
\Sm \\
\Mhm \\
\Rb
\end{array} \right).
\edm
Then $\textrm{rank}\, \Fm = n$.
\end{corollary}

\emph{Proof.}
As before, every $(2j-1)^{th}$ row of $\Fmu$ is equal to the $j^{th}$ column of $\Fm$ and every $2j^{th}$ row of $\Fmu$ is (-1) times the $(n+j)^{th}$ column of $\Fm$.  Thus, 
\bdm
\textrm{rank}\, \Fm = \textrm{rank}\, \Fmu^{t} = \textrm{rank}\, \Fmu = n
\edm
That is, $\textrm{rank}\, \Fm = n$.  \hfill $\blacksquare$

\section{How to Test for a Solution} \label{sol}
Let 
\begin{eqnarray*}
f_{i}(z_{1},\ldots,z_{n}) & = & \sum_{j=1}^{n}\Big( (p_{i}m'_{ij} + q_{i}l'_{ij}) \llg(z_{j}) + (p_{i}m''_{ij} + q_{i}l''_{ij}) \llg(1-z_{j}) \Big) \\
 & & \qquad \qquad \qquad \quad - c_{si}\pi {i}\  \qquad \quad  (i=1, \ldots, k-h) \\
f_{i}(z_{1},\ldots,z_{n}) & = & \sum_{j=1}^{n}(m'_{ij} \llg(z_{j}) + m''_{ij} \llg(1-z_{j})) - c_{mi} \pi {i} \\
 & &  \qquad\qquad\qquad\qquad \qquad\qquad\quad  (i=k-h+1, \ldots, k) \\
f_{i}(z_{1},\ldots,z_{n}) & = & \sum_{j=1}^{n}(r'_{ij} \llg(z_{j}) + r''_{ij} \llg(1-z_{j})) - c_{i} \pi {i} \\
 & & \qquad \qquad\qquad\qquad\qquad \qquad\qquad\qquad \,\,\, (i=k + 1, \ldots, n)
\end{eqnarray*}
and let 
\begin{eqnarray*}
f:\C^{n} & \to & \C^{n} \qquad \textrm{  such that}\\
z=(z_{1},\ldots, z_{n}) & \to & f(z) = (f_{1}(z),\ldots, f_{n}(z)).
\end{eqnarray*}
Then let

\bdm
\begin{array}{llll}
t'_{ij}=p_{i}m'_{ij}+q_{i}l'_{ij} & t''_{ij}=p_{i}m''_{ij}+q_{i}l''_{ij} & t'''_{i}=c_{si} & (i=1, \ldots, k-h)\\
t'_{ij}=m'_{ij} & t''_{ij}=m''_{ij} & t'''_{i}=c_{mi} & (i=k-h+1, \ldots, k)\\
t'_{ij}=r'_{ij} & t''_{ij}=r''_{ij} & t'''_{i}=c_{i} & (i=k+1, \ldots, n).
\end{array}
\edm
The resulting components of $f$ are 
\beqn
f_{i}(z_{1},\ldots,z_{n}) = \sum_{j=1}^{n}(t'_{ij} \llg(z_{j}) + t''_{ij} \llg(1-z_{j})) - t'''_{i} \pi {i}  \qquad  (i=1, \ldots, n).
\eeqn
Then $\partialfi{i}{z} = \frac{t'_{ij}}{z_{j}} - \frac{t''_{ij}}{1-z_{j}} \qquad \textrm{ for } 1\leq i\leq n,$ so \label{fijpartial}
\beqn
\partialf{z} = \Big( {\frac{t'_{1j}}{z_{j}}} - {\frac{t''_{1j}}{1-z_{j}}},\ldots,{\frac{t'_{nj}}{z_{j}}} - {\frac{t''_{nj}}{1-z_{j}}} \Big).
\eeqn
Let $H = \C_{+}^{n},$ the upper half plane in $\C^{n}.\,\,H$ is open in $\Cn$.  Each $f_{i}$ is holomorphic on $H,$ so $f$ is holomorphic on $H$~\cite{Range}. Thus $f$ is smooth on $H,$ with the derivative of $f$ at $z,\,f'(z) = \Big( \partialfi{i}{z} \Big)_{1\leq i,j\leq n}$~\cite{Whit}, being well defined on $H$.  Since we are only working with manifolds where SnapPea finds an approximate solution to $f$ in $\C_{+}^{n},$ there is an $a \in \C_{+}^{n}$ such that $f(a) = b \aand b$ is extremely close to $0 \in \Cn$.  We know that $\textrm{det} f'(a)\ne 0~\cite{Choi},$ so rank $f'(a) = n \aand f$ is regular at $a$.  Then $f'(a)^{-1}$ exists.  Let
\begin{eqnarray*}
\delta:\Cn & \to & \Cn\\
v & \to & |f'(a)^{-1} \cdot v|.
\end{eqnarray*}
Since $\delta,$ as a function of $v,$ is a continuous function on $\Cn,$ it will attain a maximum and minimum on the compact set $\{v \in \Cn : |v| = 1\}$.

\subsection{Kantorovich} \label{secKant}
The Kantorovich Theorem~\cite{Kant} provides a test for the solution of $f$.  The relevance of this theorem to the solution of $f$ was brought to our attention by Joan Birman after another test had been developed by us.  We thank her for telling us about it. The Kantorovich Theorem is usable in our situation because we can identify the quantities used, though this is not the case for all functions.  The test provides a sufficient condition for a manifold to have a complete hyperbolic structure.  Consequently, it is possible for a manifold to not satisfy the condition and still be complete hyperbolic.  
\begin{theorem}[Kantorovich] Let $U$ be an open neighborhood of a point, $a,$ in $\Cn$ and $f:U \to \Cn$ a holomorphic mapping with invertible derivative $f'(a)$ at $a.$  Let $hh = -f'(a)^{-1} f(a), \, \tilde{a} = a + hh \aand U_{0} = B_{|hh|}(\tilde{a})$, the open ball of radius $|hh|$ about $\tilde{a}$.  If $U_{0} \subset U$ and
\begin{enumerate}
\item The derivative $f'(z)$ satisfies the Lipschitz Condition on $U_{0},$ with Lipschitz Ratio, $L$
\item $|f(a)||f'(a)^{-1}|^{2} L \leq \frac{1}{2},$
\end{enumerate}
then $f(z) = 0$ has a unique solution in $U_{0}.$
\end{theorem}
The Kantorovich Theorem applied to our function, $f,$ works as follows.  Let $U = H.$  Given $a,$ an approximate solution to $f(z)=0,$ apply Newton's method to $f$ at $a$ to get an even better approximate solution, $\tilde{a}.$  That is, let $hh = -f'(a)^{-1} \cdot f(a) \aand \tilde{a} = a + hh = (a_{1} + hh_{1}, \ldots, a_{n} + hh_{n})$ so $\tilde{a}_{j} = a_{j} + hh_{j}.$  Then see if a Lipschitz Ratio, denoted by $L,$ can be identified for $z \in \Bhh$ so that $f'(z)$ satisfies the Lipschitz condition on $U_{0}$ with $L.$  One way to do this is to find an upper bound, $c_{ijk},$ on the second partials, $|\partial_{i}\partial_{j}f_{k}(z)|$ for $1 \leq i,j,k \leq n$ for $z \in \Bhh,$ and let $L = \sqrt{\sum_{1 \leq i,j,k \leq n}(c_{ijk})^{2}}~\cite{Kant}.$  This works for us, but in general, the major stumbling block to using this theorem is the difficulty in finding this $L$.  Here, $|f'(a)^{-1}|,$ the norm of $f'(a)^{-1},$ can be either the supremum norm, which we will denote by $|f'(a)^{-1}|_{\mathrm{sup}},$ or the length norm, referred to as $|f'(a)^{-1}|_{\mathrm{len}},$ where
\bdm
|f'(a)^{-1}|_{\mathrm{sup}} = \textrm{sup}_{|v| = 1} |f'(a)^{-1} \cdot v|
\edm
and if a component of $f'(a)^{-1}$ is denoted by $h_{ij},$
\bdm
|f'(a)^{-1}|_{\mathrm{len}} = \sqrt{\sum_{1 \leq i,j \leq n}|h_{ij}|^{2}}.
\edm
Now substitute values in the inequality found in the second part of the Kantorovich Theorem and see if they pass the test.  If so, there is a solution in $\Bhh.$

\subsubsection{Calculate $|f'(a)^{-1}|$}
\begin{description}
\item [Supremum Norm:$|f'(a)^{-1}|_{\mathrm{sup}}\,\,$]

Let
\bdm
B = \{f'(a)^{-1} \cdot v: |v| = 1\} = \{ w \in \Cn:|f'(a) \cdot w| = 1\}.
\edm
We look at the continuous real valued function $\mu$ on the compact set $B$ such that 
\begin{eqnarray*}
\mu:B &\to & \R \\
w &\to & |w|^{2}.
\end{eqnarray*}
Let $S = \{v\in \Cn:|v| = 1\}.$  Then $\mu$ attains a maximum at some $\tilde{w} \in B$  and the function $\delta$ will attain a maximum at some $\tilde{v} \in S$ where $\tilde{w} = f'(a)^{-1} \cdot \tilde{v}$.
Now let $A = f'(a)$.  This is a complex matrix, so 
\begin{eqnarray*}
|Aw|^{2} & = & (Aw)^{t}(\overline{Aw}) \qquad \overline{A} = \textrm{conjugate of }A \aand t = \textrm{transpose of }A\\
& = & (w^{t}A^{t})(\overline{Aw})\\
& = & w^{t}(A^{t}\overline{A})\overline{w}.
\end{eqnarray*}
Let $D = (A^{t}\overline{A})$.  This is a self adjoint matrix so it has real eigenvalues~\cite{Edw}.  Then,
\begin{eqnarray*}
B & = & \{w:|Aw| = 1\}\qquad \qquad \qquad \qquad \qquad\\
& = & \{w:|Aw|^{2} = 1\}\\
& = & \{w:w^{t} D \overline{w} = 1\}.
\end{eqnarray*}
Using the Lagrange multiplier method to maximize $\mu$ on $B$~\cite{Edw}, let
\begin{eqnarray*}
H(w_{1}, \ldots, w_{n},\lambda) & = & |w|^{2} - \lambda(w^{t} D \overline{w} -1 )\\
& = & \sum_{i=1}^{n} w_{i}\overline{w_{i}} - \lambda\Big(\sum_{i=1}^{n}w_{i}(\sum_{j=1}^{n} d_{ij}\overline{w_{j}}) - 1 \Big).
\end{eqnarray*}
In order to find a critical point for $H,$ all partials with respect to $w_{1}, \ldots, w_{n}$ and $\lambda$ must be $0$.  We set
\bdm
0 = \frac{\partial H}{\partial w_{i}} = \overline{w_{i}} - \lambda (\sum_{j=1}^{n} d_{ij}\overline{w_{j}}) \qquad (i=1, \ldots, n),
\edm
so,
\begin{eqnarray*}
0 & = & \overline{w} - \lambda D \overline{w} \\
& = & (I - \lambda D)\overline{w}\\
& = & (\frac{1}{\lambda} I - D)\overline{w}.
\end{eqnarray*}
Then $D\overline{w} = \frac{\overline{w}}{\lambda},$ making $\frac{1}{\lambda}$ an eigenvalue of $D$.  Also,
\begin{eqnarray*}
0 & = & \frac{\partial H}{\partial \lambda}\\ 
& = & w^{t}D\overline{w} - 1.
\end{eqnarray*}
Thus, $w^{t}D\overline{w} = 1,$ and substituting $\frac{\overline{w}}{\lambda}$ for $D\overline{w}$ from above, we have $w^{t} \frac{\overline{w}}{\lambda} = 1$.  That is, $w^{t} \overline{w} = \lambda$.  But $w^{t} \overline{w} = |w|^{2},$ so $\textrm{max}|w|^{2} = \textrm{max}\lambda$ such that $\frac{1}{\lambda}$ is an eigenvalue of $D$.  Then,
\begin{eqnarray*}
\textrm{max}_{\mathrm{on }B}|w|^{2} & = & \frac{1}{\textrm{smallest eigenvalue of }D}\\
& = & \frac{1}{\textrm{smallest eigenvalue of }A^{t}\overline{A}}\\
& = & \frac{1}{\textrm{smallest eigenvalue of }f'(a)^{t}\overline{f'(a)}}.
\end{eqnarray*}
By definition, $|f'(a)^{-1}|_{\mathrm{sup}} = \textrm{max}_{\mathrm{on }B} |w|,$ so
\begin{equation}
|f'(a)^{-1}|_{\mathrm{sup}} = \frac{1}{\sqrt{{\textrm{smallest eigenvalue of }f'(a)^{t}\overline{f'(a)}}}}.
\end{equation}
We calculate the eigenvalues of $f'(a)^{t}\overline{f'(a)}$ using its characteristic polynomial and take the square root of the smallest one to get $|f'(a)^{-1}|_{\mathrm{sup}}$. \label{eigen}
\item [Length Norm:$|f'(a)^{-1}|_{\mathrm{len}}\,\,$]

Let the components of $f'(a)^{-1}$ be $(h_{ij})_{1 \leq i,j \leq n}.$  Then
\beqn
|f'(a)^{-1}|_{\mathrm{len}} = \sqrt{\sum_{1 \leq i,j \leq n}|h_{ij}|^{2}}.
\eeqn
\end{description}
\subsubsection{Calculate $c_{ijk}$}
Let $z \in \Bhh.$  Then $|z - \tilde{a}| < |hh|,$ so $|z_{j} - \tilde{a}_{j}| < |hh|,$ where $z_{j} - \tilde{a}_{j} = z_{j} - ({a}_{j} + hh_{j})$ since $\tilde{a}_{j} = {a}_{j} + hh_{j}.$  Figure~\ref{kantpic} shows the situation for each $j$.  There are three tests that need to be performed before we test for the inequality in the Kantorovich Theorem.  The entire process stops and Kantorovich tells us nothing about a manifold when any of these tests fail.
\begin{figure}
\begin{picture}(250,135)(-120,-25)
\put(-26,0){\line(1,0){200}}
\put(0,-26){\line(0,1){125}}
\put(2,-9){$0$}
\put(130,-9){$1$}
\put(84,58){\textbf{.}}
\put(77,54){$a_{j}$}
\put(97,74){\circle{39}}
\put(97,74){\textbf{.}}
\put(102,72){$\tilde{a}_{j}$}
\put(90,80){\textbf{.}}
\put(90,85){$z_{j}$}
\end{picture}
\caption{Disc of radius $|hh_{j}|$ about $\tilde{a}_{j}$}  \label{kantpic}
\end{figure}
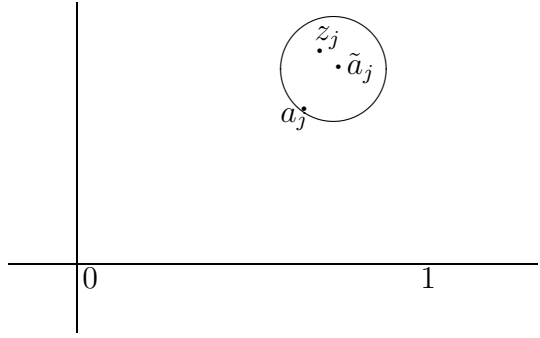 
\begin{description}
\item [Test 1]
We want a solution in $H,$ so \textbf{we require that} $\mathbf{\mathrm{Im}(\tilde{a}_{j}) > |hh|}.$  Otherwise, there are $ z \in \Bhh$ that have $\textrm{Im}(z_{j}) \leq 0,$ and the solution could be one of these $z.$
\item [Test 2]
\bdm
|(z_{j} - a_{j}) - hh_{j}| = |z_{j} - (a_{j} + hh_{j})| = |z_{j} - \tilde{a}_{j}| < |hh|.
\edm
Using triangle inequalities,
\bdm
|z_{j} - a_{j}| - |hh_{j}| \leq |(z_{j} - a_{j}) - hh_{j}|.
\edm
Therefore,$|z_{j} - a_{j}| - |hh_{j}| < |hh|,$ giving $|z_{j} - a_{j}| < |hh_{j}| + |hh|.$  But $|hh_{j}| \leq |hh|,$ so $|z_{j} - a_{j}| < 2|hh|.$  Now
\bdm
|z_{j}| = |a_{j} + (z_{j} - a_{j})| \geq |a_{j}| - |(z_{j} - a_{j})|.
\edm 
Thus,
\beqn
|z_{j}| > |a_{j}| - 2|hh|.
\eeqn
We need $|a_{j}| - 2|hh| >0$ in order to define $L,$ so \textbf{the second test is to check that} $\mathbf{ |hh| < \frac{1}{2}|a_{j}|}.$
Then,
\begin{equation}
\frac{1}{|z_{j}|} < \frac{1}{|a_{j}| - 2|hh|}.    \label{test2}
\end{equation}
\item [Test 3]
We do a similar process as in the previous test.  We already know that $|z_{j} - a_{j}| < 2|hh|.$  Only now, we use $1 - z_{j}$ instead of $z_{j},$ so
\bdm
|1 - z_{j}| = |1 - a_{j} - (z_{j} - a_{j})| \geq |1 - a_{j}| - |(z_{j} - a_{j})|.
\edm 
Hence,
\beqn
|1 - z_{j}| > |1 - a_{j}| - 2|hh|.
\eeqn
We need $|1 - a_{j}| - 2|hh| >0;$\textbf{ a third test is to see if} $\mathbf{ |hh| < \frac{1}{2}|1 - a_{j}|}.$
Then,
\begin{equation}
\frac{1}{|1 - z_{j}|} < \frac{1}{|1 - a_{j}| - 2|hh|}.   \label{test3}
\end{equation}
\item[Remainder of Calculation]
\end{description}

We are now ready to look at the second partials.  By page~\pageref{fijpartial} we see that for $z \in \Bhh,$
\bdm
\partial_{j}f_{i}(z) = \partialfi{i}{z} = \frac{t'_{ij}}{z_{j}} - \frac{t''_{ij}}{1-z_{j}}.
\edm
Therefore,
\begin{eqnarray}
\partial_{k}\partial_{j}f_{i}(z) & = & 0\qquad \textrm{for } k \neq j \\
& = & -\frac{t'_{ij}}{z_{j}^{2}} - \frac{t''_{ij}}{(1-z_{j})^{2}} \qquad \textrm{for }k = j.
\end{eqnarray}
Consequently,
\bdm
|\partial_{j}\partial_{j}f_{i}(z)| \leq \frac{|t'_{ij}|}{|z_{j}|^{2}} + \frac{|t''_{ij}|}{|1-z_{j}|^{2}}.
\edm
Combining this with Equations~\ref{test2} and~\ref{test3} yields
\beqn
|\partial_{j}\partial_{j}f_{i}(z)| \leq \frac{|t'_{ij}|}{(|a_{j}| - 2|hh|)^{2}} + \frac{|t''_{ij}|}{(|1 - a_{j}| - 2|hh|)^{2}}.
\eeqn
Using this, $c_{ijk}$ for $1 \leq i,j,k \leq n$ is defined as
\begin{eqnarray}
c_{ijk} & = & 0 \qquad \textrm{for } j \neq k  \\
c_{ijj} & = & \frac{|t'_{ij}|}{(|a_{j}| - 2|hh|)^{2}} + \frac{|t''_{ij}|}{(|1 - a_{j}| - 2|hh|)^{2}}.
\end{eqnarray}
The Lipschitz Ratio, $L,$ can now be identified as
\begin{eqnarray}
L & = & \sqrt{\sum_{1 \leq i,j \leq n}(c_{ijj})^{2}} \nonumber \\
L & = &   \sqrt{\sum_{1 \leq i,j \leq n}\Bigg(\frac{|t'_{ij}|}{(|a_{j}| - 2|hh|)^{2}} + \frac{|t''_{ij}|}{(|1 - a_{j}| - 2|hh|)^{2}}\Bigg)^{2}}.
\end{eqnarray}

The theorem can finally be applied, testing to see if $|f(a)||f'(a)^{-1}|^{2}L \leq \frac{1}{2}.$  Since $b = f(a),$ this can be rewritten as
\bdm
|b| \leq \frac{1}{{2}|f'(a)^{-1}|^{2}{L}}.
\edm
We really have two tests, one using the supremum norm and the other using the length norm.  This completes the last part of the proof of Theorem~\ref{th}.

\section{Examples} \label{ex}
The methods presented are implemented by the use of two programs: SNAP to get information about the manifold and Pari-Gp~\cite{PGP} to do calculations.  We use Pari-Gp instead of Mathematica because of its high level of precision.  Page~\pageref{temp} begins the template for a program that is written in an edit file and then copied into Pari-Gp for execution.  The template needs to be adjusted for information gotten from SNAP.  Assume we have a manifold file in SNAP for the manifold, $M.$  Once SNAP is open, read in the file for processing.  The ``pr sol" command will print the type of solution SNAP has found.  A geometric solution means that the solution is in $H.$  Any other response is useless here, so there is no need to go any further.  Assuming it is geometric, proceed with setting up the template.  Issue the ``pr sh" command.  SNAP will return the transpose of a vector representing an approximate solution to our set of $n$ equations for $M.$  The number of components of the vector will be equal to $n,$ the number of tetrahedra in the triangulation.  Then copy this vector from SNAP to the template, replacing $[a_{1}, \ldots, a_{n}],$ so that $a$ now has the value of our approximate solution.  The tilde at the end of the Snap response must be eliminated so that $a$ appears as a $1 \times n$ matrix.  Next comes the ``pr fill" command.  SNAP will display a $(n + k) \times (2n +1)$ matrix where the components of each row are the coefficients of a cusp or consistency equation.  Assuming $M$ is the result of Dehn filling on $h$ out of $k$ cusps, the first $k - h$ rows represent the cusp surgery equations, the next $h$ rows are the meridian completeness equations for the unsurgered cusps, and the last $n$ rows are all the consistency equations before any have been eliminated.  If all cusps are unsurgered, $h = k,$ so the first $k$ rows are all meridianal completeness equations.  This command to print filling equations is closely related to the ``pr gl" command which prints the gluing equations.  This latter display presents the $k$ meridianal followed by the $k$ longitudinal completeness equations for the original $k$ cusps before any surgery, and then all of the $n$ consistency equations.  But it is simpler to use the filling equations, even for manifolds where no surgery has been done.  Copy this matrix from Snap to the template, initializing the matrix $FG$.  The script will then create the matrices $F \aand G,$ where $F$ consists of the first $k$ rows of $FG$ and $G$ consists of the last $n$ rows of $FG.$  The rows of $F$ are linearly independent and the program selects $n - k$ rows from $G$ so that when added to $F,$ the resulting matrix has rank $n.$ 

The only further adjustments may be Pari-Gp punctuation to reflect line continuation.  In order to tell Pari-Gp to ignore an end of line from the text editor, a ``$\setminus$" followed immediately by Return must end that text line.  This is needed with a large vector or matrix, so it will probably be needed once the values for $a,\,F \aand G$ are copied into the template.
\subsection{Template}
/* set precision to $60$ (or higher for very large manifolds) from the default of */ \\
$\setminus$p $60$\\   \label{temp}
/* read the file FILENAME into SNAP */ \\
/* see that there are h unsurgered cusps */ \\
/* (,) (,) (,) (,) */ \\
/* print shapes - the triangulation has $n$ tetrahedra */\\
/* enter the shapes as a vector, so it is regarded as a $1\times n$ matrix */\\
a $=$ \\
$\bbb \bbb [a_{1},\ldots, a_{n}]$ \\
/* find $n,$ the number of tetrahedra  */\\
n $=$ matsize(a)[2]\\
/* print  filling equations and use this to initialize the matrix FG.  The first  \eol\\
$\bbb k-h$ equations are cusp surgery equations, followed by the $h$ meridianal \eol\\
$\bbb$ completeness equations, and finally all of the $n$ consistency equations. */\\
FG $=$ \\
 $\bbb \bbb [x_{11}, \ldots, x_{1(2n+1)}; \ldots; x_{(n+k)1}, \ldots, x_{(n+k)(2n+1)}]$\\
/* find $n+k,$ the number of equations derived from the Snap command \eol \\
$\bbb$ "pr fill" */\\
numalleq $=$ matsize(FG)[1]\\
/* find total number of cusps, $k$ */\\
k $=$ numalleq - n\\
/* initialize F, the cusp equations matrix, using the first $k$ equations from \eol \\
$\bbb$ FG */\\
F $=$ matrix(k,2*(n) +1,ii,jj,FG[ii,jj])\\
/* initialize G, the matrix of all consistency equations, using the last $n$ \eol\\
$\bbb$ equations from FG */ \\
G $=$ matrix(n,2*(n) +1,ii,jj,FG[k+ii,jj])\\
/* define matrix $H$ by eliminating the last column of $F$ representing the  \eol \\
$\bbb \pi i$ coefficient */ \\
H $=$ matrix (k,2*n,i,j,F[i,j])\\
/* define matrix $K$ by eliminating the last column of $G$ representing the  \eol \\
$\bbb \pi i$  coefficient */ \\
K $=$ matrix (n,2*n,i,j,G[i,j])\\
/* redefine $F$ and $H$ by adding rows to them from $G$ and $K$ respectively \eol \\
$\bbb$ until the rank of $F$ and $H$ are both $n$ */ \\
r $=$ 1 \\
v(r) $=$ vector( (2*n)$+$1,l,G[r,l] ) \\
t(r) $=$ vector( 2*n,l,K[r,l] )  \\
while( n $-$ matrank(H)$ \bb \& \& \bb$(n$+$1$-$r), if( (matrank(concat(F,v(r))) $\bs \eolm \\
$\BbB -$ matrank(F))$ \bb \& \& \bb$(matrank(concat(H,t(r))) $-$ matrank(H)), \eol \\
$\BbB $ (F $=$ concat(F,v(r)))$ \bb \& \& \bb$ (H $=$ concat(H,t(r))), r$=$r$+$1)) \\
eval(F) \\
eval(H) \\
/* set up the filling equations as log functions evaluated at $a$ */\\
f(i) $=$ sum( j $=$ 1, n, F[i,j]*log(a[j]) ) $+$ ( j = 1, n, F[i,n+j]*log(1-a[j]) ) \eol \\
$\BBb \BBb +$ F[i,(2*n)$+$1]*Pi*I \\
/* define the vector $b$ in $\Cn$ */ \\
b $=$ vector( n, i, f(i) ) \\ 
/* identify the norm of $b$ */ \\
normb $=$ sqrt( norml2(b) ) \\
/* identify $A,$ the derivative matrix for $f$ at $a$ */\\ 
g(i,j) $=$ ( F[i,j]/a[j] ) $-$ ( F[i,n+j]/(1-a[j]) ) \\
A $=$ matrix( n, n, i, j, g(i,j) ) \\
/* check that determinant of $A$ is not zero */\\
matdet(A) \\

/* KANTOROVICH PROCESSING */ \\ 
/* find eigenvalues for $D=(\textrm{transpose of }A) * (\textrm{conjugate of }A)$ */ \\
D $=$ mattranspose(A)*conj(A) \\
wapprox $=$ polroots( charpoly(D,x) ) \\
w $=$ real( wapprox ) \\ 
/* change $b$ into a matrix to do matrix multiplication */ \\
B $=$ matrix(n,1,j,i,b[j]) \\
/* define the vector $hh$ and find its length, normhh */ \\
hhh $=$ $-$(A)\verb+^+($-1$)*(B) \\
hh $=$ vector(n, j, hhh[j,1]) \\
normhh $=$ sqrt(norml2(hh)) \\
/* perform the first three tests to see if this method is applicable */ \\
atilde = a $+$ hh \\
/* test 1 to see if fat solution; if j $>$ n */ \\
for (j $=$ 1, n, if(normhh $<$ imag(atilde[j]), , \eol \\
$\BBb \BBb$ error(''failure at atilde['', j, '']''))) \\
/* test 2 to see if $c_{ijj}$ can be defined */ \\
for(j $=$ 1, n, if(normhh $<$ (1/2)*abs(a[j]), , \eol \\
$\BBb \BBb$ error(''failure at atilde['', j, '']''))) \\
/* test 3; other test to see if $c_{ijj}$ can be defined */ \\
for(j $=$ 1, n, if(normhh $<$ (1/2)*abs(1 $-$ a[j]), , \eol \\
$\BBb \BBb$ error(''failure at atilde['', j, '']''))) \\%
/* identify the Lipschitz ratio, Lips */ \\
c(i,j) $=$ (abs(F[i,j])/(abs(a[j]) $-$ 2*normhh)\verb+^+2) \eol \\
$\BBb \BbB +$ (abs(F[i,j$+$n])/(abs(1$-$a[j]) $-$ 2*normhh)\verb+^+2) \\ 
Lips $=$ sqrt( sum( j $=$ 1, n, sum(i $=$ 1, n, c(i,j)\verb+^+2) ) ) \\
/* identify normAinv, the norm of A\verb+^+($-1$), using the definition of matrix \eol \\
$\bbb$ norm as the supremum of A\verb+^+($-1$)$v$ for $v$ on the $n$-sphere */ \\
normAinv $=$ 1/sqrt( vecmin(w) ) \\ 
/* do the Kantorovich tests */ \\
/* find the value that the norm of $b$ must be less than or equal to with  \eol \\
$\bbb$ respect to the supremum norm */ \\
1/(2 * (normAinv)\verb+^+2 * Lips) \\
normb $<=$ 1/(2 * (normAinv)\verb+^+2 * Lips) \\
/* find the length norm and the value that the norm of $b$ must be less than\eol \\
$\bbb$ or equal to with respect to the length norm */ \\
sqrt(norml2(A\verb+^+($-1$))) \\
1/(2 * norml2(A\verb+^+($-1$)) * Lips) \\
normb $<=$ 1/(2 * norml2(A\verb+^+($-1$)) * Lips) \hfill $\blacksquare$

\subsection{Using the Template}
The template is ready to be used.  If you want a copy of what has happened, first turn the log on in Pari-Gp by typing ``$\setminus$l logfilename".  Then copy the adjusted template to Pari-Gp, wait for the run to complete, and open the log file to see the results.  Make sure that there are no error messages from the qualification tests described.  If there are, any further results are of no value. If there are no error messages, a response of ``$1$" to either of the Kantorovich inequalities will indicate the manifold is complete hyperbolic.  A copy of the template can be found at~\cite{HM}.  We now look at three examples.  Each example will have two sets of data.  The first comes from SNAP and the second is the result of calculations in Pari-Gp.  The vectors and matrices are printed as they appear in SNAP.  When one of them extends beyond one line, it is edited once copied into the template to add the line continuation character, ``$\setminus,$" after each line before its end.  The Pari-Gp data has been shortened to 40 decimal places from the calculated precision of 60 decimal places so as to fit on one line since in these examples, it has no effect on understanding the results.
\begin{enumerate}
\item \underline{\emph{FIGURE $8$ KNOT COMPLEMENT}}\\  
The simplest is the figure $8$ knot complement.  We know~\cite{Thur} that this is complete hyperbolic already.  However, only sufficiency conditions have been presented here, so it is nice to see that a manifold we know to be complete hyperbolic does not fail the test.
\\[.14in]
QUANTITIES FROM SNAP\\
$n = 2$\\
$h = k = 1$ \\
$a = $\\
{[}0.5000000000000000000000000000+0.8660254037844386467637231707*I,\\
0.5000000000000000000000000000+0.8660254037844386467637231707*I]\\
$F = $[1, 0, 0, 1, 0]\\
$G = $[2, -1, -1, 2, 0; -2, 1, 1, -2, 0] \\[.14in]
Pari-Gp CALCULATIONS\\
$|b| = 1.296666384352891444530724934775173278518 E-28$\\
$L = 4.472135954999579392818347339211785668123$ \\
$|f'(a)^{-1}|_{\mathrm{sup}} = 1.592226038754547070932399593119376104348$ \\
$|f'(a)^{-1}|_{\mathrm{len}} = 1.632993161855452065464856049716587347937$ \\
$\frac{1}{2 {|f'(a)^{-1}|_{\mathrm{sup}}^{2}} L } = 0.04410070808503045666350407221846082500302$\\
$\frac{1}{2 {|f'(a)^{-1}|_{\mathrm{len}}^{2}} L } = 0.04192627457812105680767200627679720162466$  

\item \underline{\emph{$(9872,11111)$ DEHN SURGERY: WHITEHEAD LINK COMPLEMENT}}\\  
The Whitehead link complement is known to be complete hyperbolic~\cite{N-R}.  This example considers Dehn surgery on only one of the two cusps of the Whitehead link complement.
\\[.14in]
QUANTITIES FROM SNAP\\
$n = 4$\\
$h = 1 \aand k = 2$ \\
$a = $\\
{[}0.9999343700073827649570992430+1.000170536257729817727630077*I,\\
0.4999147436597508540443693049+0.4999671844066970777583211769*I,\\
0.5000852675298210651958243937+0.5000328032070212542658981140*I,\\
0.4999147436597508540443693049+0.4999671844066970777583211769*I]\\
$F = $\\
{[}20983, 0, -9872, 0, -9872, 11111, -1239, 20983, -2;\\
 0, 0, 0, 1, 1, -1, 0, 0, 0]\\
$G = $\\
{[}1, 1, 1, 1, 1, -2, 0, 0, -1; 0, -1, -1, -1, -1, 1, 1, 1, 1;\\
-1, 1, 1, 1, 1, 0, -2, -2, -1; 0, -1, -1, -1, -1, 1, 1, 1, 1]\\[.14in]
Pari-Gp CALCULATIONS\\
$|b| = 6.290546043622649509854067366063508951285 E-24$\\
$L = 56237.01131396100111291495604741250466464$ \\
$|f'(a)^{-1}|_{\mathrm{sup}} = 1.063909899076773471157618529051471308315$ \\
$|f'(a)^{-1}|_{\mathrm{len}} = 1.235415661324873497175222236812823735348$ \\
$\frac{1}{2 {|f'(a)^{-1}|_{\mathrm{sup}}^{2}} L } = 0.000007854853193291278165225494981053686965848$\\
$\frac{1}{2 {|f'(a)^{-1}|_{\mathrm{len}}^{2}} L } = 0.000005825343870778317976532920417278552662252$

\item \underline{\emph{LARGELINK COMPLEMENT}}\\  
This is the smaller of two extremely large link complements.  See figure~\ref{small}.  It has $32$ tetrahedra and $4$ cusps.  The other one has 57 tetrahedra and 11 cusps.  These two links are used by Leininger~\cite{Lein} to construct other knots and links by cut and paste methods, and then looking at their covers.  For any even integer $g>0,$ we eventually get from Largelink a two component link whose complement in $S^{3}$ contains an embedded totally geodesic surface of genus $g$.  The importance of Largelink is that prior to this, such embedded surfaces could only be found in the complement of links with more than two components.
 \\[.14in]
 \begin{figure} 
 \begin{center} 
   \includegraphics[height=1.8in]{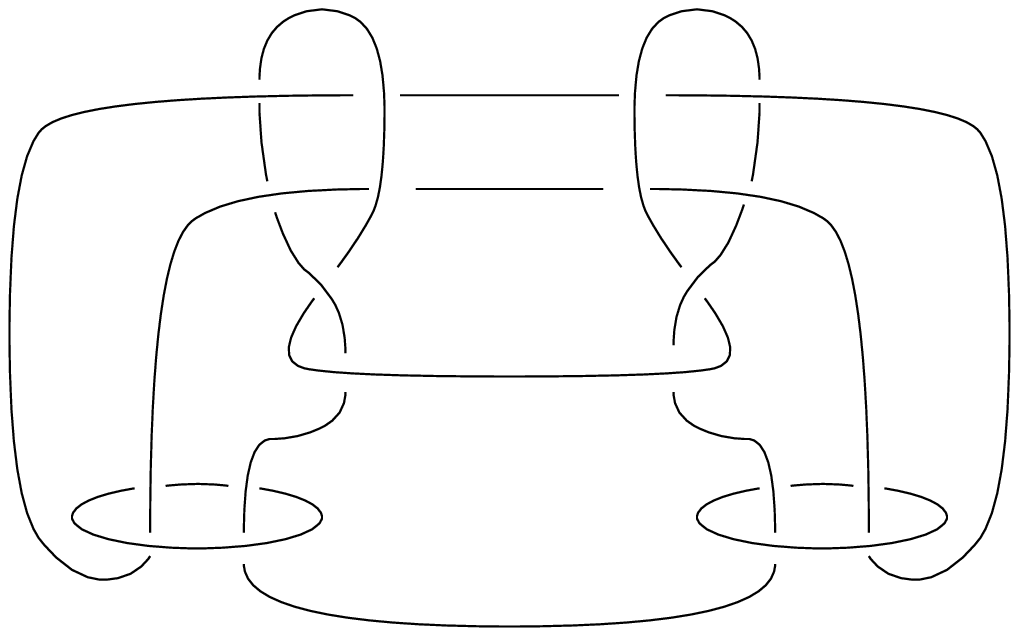}
 \caption{The Link Largelink} \label{small} 
 \end{center} 
 \end{figure} 
QUANTITIES FROM SNAP\\
$n = 32$\\
$h = k = 4$ \\
$a = $\\
{[}5.431680776271168985E-77+1.043190149785894973378994944*I,\\
0.4788708557877967957032308372+0.4995533597773714501527030266*I,\\
-4.471822153042346518E-77+0.9585980084313877504633692171*I,\\
0.5211291442122032042967691627+0.4995533597773714501527030266*I,\\
0.2929970420861826752219808548+1.473911044296957810855392169*I,\\
-0.4509782171525463654321193064+1.200765444220459728291241593*I,\\
1.000000000000000000000000000+0.9300613056344272435239940348*I,\\
0.4638110047891777790136229363+0.4986886369525889712130902195*I,\\
0.2371128008078259554449875702+0.6313317290357266968810549581*I,\\
0.3060049499572359024927254903+0.4485055715132523515850447879*I,\\
0.8585412796265143611585133046+1.027932770073775455116265474*I,\\
0.5000000000000000000000000000+0.4792990042156938752316846085*I,\\
0.4375240155821198504166057790+0.8813536566549109733830907053*I,\\
1.042258288424406408593538325+0.9991067195547429003054060533*I,\\
0.6696776343174901312972923995+0.7426519144895069642064793083*I,\\
0.4518888703362351094400929330+0.9102903934144040876554144906*I,\\
0.03927009472823897821842546946+1.571359648665194056162768058*I,\\
-8.46343996E-78+0.9585980084313877504633692171*I,\\
1.000000000000000000000000000+1.043190149785894973378994944*I,\\
0.02451089142372681728394034675+0.5982980953722294364585245396*I,\\
0.9621628947892310086730291057+0.3083453854492406606721071067*I,\\
0.7354295168083648566686302069+0.5515583107626382072967381105*I,\\
0.6213864977872760582396031709+0.4161571993484503024065288682*I,\\
0.6756917822944407062548472825+0.1978399260627268524593119332*I,\\
0.5213536432299720005859050458+1.346701507985612627940863123*I,\\
0.2659365860052524158474189000+0.5690611275237113909012011313*I,\\
0.8916797222785394793793396465+0.5330292860478110834980601119*I,\\
0.4489838724616496515202858898+0.4713823217067450930880172825*I,\\
0.5364433482241135276673629307+0.6234802797569418514001720639*I,\\
-0.4489884234808609710528328543+0.3884305318039174460267193001*I,\\
1.000000000000000000000000000+0.9585980084313877504633692171*I,\\
-8.49494342E-77+1.043190149785894973378994944*I]

$F = $\\
{\tiny {[}0, 0, 0, 0, 0, 0, 0, 0, 0, 0, 0, -1, 0, 0, 0, 0, 0, 1, 0, 0, 0, 0, 0, 0, 0, 0, 0, 0, 0, 0, 0, 0, 0, 0, 0, 0, 0, 0, 0, 0, 0, 0, 0, 1, 0, 0, 0, 0, 0, -1, 0, 0, 0, 0, 0, 0, 0, 0, 0, 0, 0, 0, 0, 1, 0; 
0, 0, 0, 0, 0, 0, 0, 0, 0, 0, 0, 0, 0, 0, 0, 0, 0, -1, -1, 0, 0, 0, 0, 0, 0, 0, 0, 0, 0, 0, 0, 0, 0, 0, 0, 0, 0, 0, 0, 0, 0, 0, 0, 0, 0, 0, 0, 0, 0, 1, 0, 0, 0, 0, 0, 0, 0, 0, 0, 0, 0, 0, 0, 0, 1; 
0, 1, 0, 0, 0, 0, 0, 0, 0, 0, 0, 0, 0, 0, 0, 0, 0, 0, 0, 0, 0, 0, 0, 0, 0, 0, 0, 0, 0, 0, 0, 0, 1, 0, 0, 0, 0, 0, 0, 0, 0, 0, 0, 0, 0, 0, 0, 0, 0, 0, 0, 0, 0, 0, 0, 0, 0, 0, 0, 0, 0, 0, 0, 0, 0; 
0, 0, 0, 0, 0, 0, 0, 0, 1, 0, 0, 0, 0, 0, -1, 0, -1, -1, 1, -1, -1, 0, 1, 1, 1, 2, 0, 1, 0, -1, -1, 0, 0, 0, 0, 0, 0, 0, 0, 1, 1, 0, 0, 0, 0, 0, 0, 0, 1, 1, 0, 1, 0, -2, 0, 0, -1, 0, 0, 0, 0, 1, 1, -1, 1] \par}

$G = $\\
{\tiny {[}0, 1, -1, 0, 0, 0, 0, 0, 0, 0, 0, 0, 0, 0, 0, 0, 0, 0, 0, 0, 0, 0, 0, 0, 0, 0, 0, 0, 0, 0, 0, 0, 1, 0, 1, 1, 0, 0, 0, 0, 0, 0, 0, 0, 0, 0, 0, 0, 0, 0, 0, 0, 0, 0, 0, 0, 0, 0, 0, 0, 0, 0, 0, 0, 1; 
0, 0, 0, 0, 0, 0, 0, 0, 0, 0, 0, 1, 0, 0, 0, 0, 0, -1, 0, 0, 0, 0, 0, 0, 0, 0, 0, 0, 0, 0, 0, 0, -1, 0, 0, -1, 0, 0, 0, 0, 0, 0, 0, -1, 0, 1, 0, 0, 0, 1, 0, 0, 0, 0, 0, 0, 0, 0, 0, 0, 0, 0, 0, -1, 0; 
-1, 0, 0, 1, -1, 0, 0, 0, 0, 0, 0, 0, 0, 0, 0, 0, 0, 0, 0, 0, 0, 0, 0, 0, 0, 0, 0, 0, 0, 0, 0, 0, 1, -1, -1, 0, 1, 0, 0, 0, 0, 0, 0, 0, 0, -1, 0, 0, 0, 0, 0, 0, 0, 0, 0, 0, 0, 0, 0, -1, 0, -1, 0, 0, 0; 
0, -1, 1, 0, 0, 0, 0, 0, 0, 0, 0, 0, 0, 0, 0, 0, 0, 0, 0, 0, 0, 0, 0, 0, 0, 0, 0, 0, 0, 0, 0, 0, -1, 0, -1, -1, 0, 0, 0, 0, 0, 0, 0, 0, 0, 0, 0, 0, 0, 0, 0, 0, 0, 0, 0, 0, 0, 0, 0, 0, 0, 0, 0, 0, -1; 
1, -1, 0, 0, 0, 0, 0, 0, 0, 0, 0, 0, 0, 1, 0, 0, 0, 0, 0, 0, 0, 0, 0, 0, 0, 0, 1, 1, 0, 0, 0, 0, 0, 1, 0, 0, 0, 0, 0, 0, 0, 0, 0, 0, 0, 0, 0, -1, 0, 0, 0, 0, 0, 0, 0, 0, 0, 0, 0, 0, 0, 0, 0, 0, -1; 
0, 1, -1, 0, 0, -1, 0, 0, 0, 0, 0, 0, 0, 0, 0, 0, 0, 0, 0, 0, 0, 0, 0, 0, -1, 0, 0, -1, 0, 0, 0, 0, 0, 0, 1, 0, -1, 1, 0, 0, -1, 0, 0, 0, 0, 0, 0, 0, 0, 0, 0, 0, 0, 0, 0, 0, 1, 0, 0, 1, 0, 0, 0, 0, 2; 
0, 0, 1, -1, 1, 0, 0, 0, 0, 0, 0, 0, 0, -1, 0, -1, 0, 0, 0, 0, 0, 0, 0, 0, 0, 0, 0, 0, 0, 0, 0, 0, 0, 0, 0, 1, 0, 0, 0, 0, 0, 0, 0, 0, 0, 1, 0, 1, 0, 0, 0, 0, 0, 0, 0, 0, 0, 0, 0, 0, 0, 0, 0, 0, 1; 
0, 0, 0, 0, 0, 0, 0, 0, 0, 0, 0, 0, 0, 0, -1, 0, 0, 0, 0, 0, 0, 0, 0, 0, 0, 0, 0, 0, 0, -1, 0, 0, 0, 0, 0, 0, -1, 0, 0, 0, 0, -1, 0, 0, 0, 0, 1, -1, -1, 0, 0, 0, 0, 0, 0, 0, 0, 0, 0, 0, -1, 1, 0, 0, 0; 
0, 0, 0, 0, -1, 0, 0, 0, 0, 0, 0, 0, 0, 0, 0, 1, 0, 0, 0, 0, 0, 0, 0, 1, 0, 0, 0, 0, 0, 0, 0, 0, 0, 0, 0, 0, 1, -1, 0, -1, 0, 0, -1, 0, 0, 0, 0, 0, 0, 0, 0, 0, 0, 0, 0, 0, 0, 0, -1, 0, 0, 0, 0, 0, -1; 
0, 0, 0, 0, 1, 1, 0, 0, 0, 0, 0, 0, 0, 0, 0, 0, 0, 0, 0, 0, 0, 0, 0, 0, 0, 0, 0, 0, 0, 1, 0, 0, 0, 0, 0, 0, 0, 0, 0, 0, 0, 0, 0, 0, 0, 0, 0, 0, 0, 0, 0, -1, 0, 0, 0, 0, 0, 0, 0, 0, 0, 0, 0, 0, -2; 
0, 0, 0, 0, 0, -1, -1, 1, -1, 0, 0, 0, 0, 0, 0, 0, 0, 0, 0, -1, 0, 0, 1, 0, 0, 0, 0, 0, 0, 0, 0, 0, 0, 0, 0, 0, 0, 1, 1, 0, 1, 0, 0, 0, 0, 0, 0, 0, 0, 0, 0, 1, 0, 0, 0, 0, -1, 0, 0, 0, 0, 0, 0, 0, 2; 
0, 0, 0, 0, 0, 0, 0, 0, -1, 0, 0, 0, 0, 0, 0, 0, 0, 0, 0, 1, 0, -1, 0, 0, 0, 0, 0, 0, 0, 0, 0, 0, 0, 0, 0, 0, 0, -1, 0, 0, 1, 0, 0, 0, 0, 0, 0, 0, 0, 0, 0, 0, -1, 1, 0, 0, 0, 0, 0, 0, 0, 0, 0, 0, 0; 
0, 0, 0, 0, 0, 1, 0, -1, 1, 0, 0, 0, 0, 0, 0, 0, 0, 0, 0, 0, 0, 0, 0, 0, 0, 0, 0, 0, 0, 0, 0, 0, 0, 0, 0, 0, 0, 0, -1, 1, 0, 0, 0, 0, 0, 0, 0, 0, 0, 0, 0, 0, 0, 0, 0, 0, 0, 0, 0, 0, 0, 0, 0, 0, -1; 
0, 0, 0, 0, 0, 0, 1, 0, 0, 0, 0, 0, 1, 0, 0, 0, 0, 0, 0, 0, 0, 1, -1, 0, -1, 0, 0, 0, 0, 0, 0, 0, 0, 0, 0, 0, 0, 0, 0, 0, -1, 0, 0, 0, 0, 0, 0, 0, 0, 0, 0, 0, 0, 0, 1, 0, 1, -1, 0, 0, 0, 0, 0, 0, 0; 
0, 0, 0, 0, 0, 0, -1, 0, 0, 0, 0, 0, 0, 0, 0, 0, 0, 0, 0, 0, 1, 0, 1, 1, 0, 0, 0, 0, 0, 0, 0, 0, 0, 0, 0, 0, 0, 0, 0, -1, 0, -1, 0, 0, 0, 0, 0, 0, 0, 0, 0, -1, 0, 0, -1, 0, 0, 0, 0, 0, 0, 0, 0, 0, -1; 
0, 0, 0, 0, 0, 0, 1, 0, 0, 0, 0, 0, 0, 0, 0, 0, -1, 0, 0, 0, 0, 0, 0, 0, 0, 0, 0, 0, 1, 0, 0, 0, 0, 0, 0, 0, 0, 0, 0, 1, 0, 0, 1, 0, 0, 0, 0, 0, 1, 0, 0, 0, 0, 0, 0, 0, 0, 0, 0, 0, 0, 0, 0, 0, 1; 
0, 0, 0, 0, 0, 0, 0, 0, 1, 0, 0, 0, 0, 0, 0, 0, 0, 0, 1, 0, 0, 0, 0, 0, 1, 1, 0, 0, 0, 0, -1, 0, 0, 0, 0, 0, 0, 0, 0, 0, 0, 0, 0, 0, 0, 0, 0, 0, 0, 0, 0, 0, 0, -1, 0, 0, 0, 0, 0, 0, 0, 0, 1, 0, -1; 
0, 0, 0, 0, 0, 0, 0, 0, 0, 1, 0, 0, 0, 0, 1, 0, 0, 0, 0, 0, 0, 0, 0, 0, 0, 0, 1, 0, 0, 0, 0, 0, 0, 0, 0, 0, 0, 0, 0, 0, 0, 0, -1, 0, 0, 0, 0, 0, -1, 0, 0, 0, 0, 0, 0, -1, 0, 0, 0, 0, -1, 0, 0, 0, -2; 
0, 0, 0, 0, 0, 0, 0, 0, 0, -1, -1, 0, 0, 0, 0, 0, 1, 0, 0, 0, 0, 0, 0, -1, 0, 0, 0, 0, 0, 0, 0, 0, 0, 0, 0, 0, 0, 0, 0, 0, 0, 1, 1, 0, 0, 0, 0, 0, 0, 0, 0, 0, 0, 0, 0, 1, 0, 0, 0, 0, 0, 0, 0, 0, 1; 
0, 0, 0, 0, 0, 0, 0, 0, 0, -1, 0, 0, -1, 0, 0, 0, 0, 0, 0, 0, -1, 0, 0, 0, 0, 0, 0, 0, 0, 0, 0, 0, 0, 0, 0, 0, 0, 0, 0, 0, 0, 1, 0, 0, 1, 0, -1, 0, 0, 0, 0, 0, 1, -1, 0, 0, 0, 0, 0, 0, 0, 0, 0, 0, 1; 
0, 0, 0, 0, 0, 0, 0, 0, 0, 1, 0, 0, 0, 0, -1, 0, -1, 0, 0, 0, 0, 0, 0, 0, 0, 1, 0, 0, 0, 0, 0, 0, 0, 0, 0, 0, 0, 0, 0, 0, 0, 0, 0, 0, -1, 0, 1, 0, 1, 0, 0, 0, -1, 0, 0, 0, 0, 0, 0, 0, 0, 0, 0, 0, 0; 
0, 0, 0, 0, 0, 0, 0, 0, 0, 0, 1, 0, 0, 0, 0, 0, 0, 0, 0, 1, 0, 0, -1, 0, 0, 0, 0, 0, -1, 0, 0, 0, 0, 0, 0, 0, 0, 0, 0, 0, 0, 0, 0, 0, 0, 0, 0, 0, 0, 0, 0, 0, 0, 0, 1, -1, 0, 0, 0, 0, 1, -1, 0, 0, 0; 
0, 0, 0, 0, 0, 0, 0, 0, 0, 0, 0, -1, 1, 0, 1, 1, 0, 0, 0, 0, 0, 0, 0, 0, 0, 0, 0, 0, 0, 0, 0, 0, 0, 0, 0, 0, 0, 0, 0, 0, 0, 0, 0, 1, 0, -1, 0, 0, 0, 0, 0, 0, 0, 0, 0, 0, 0, 0, 0, 0, 0, 0, 0, 0, -1; 
0, 0, 0, 0, 0, 0, 0, 0, 0, 0, 0, -1, 0, 0, 0, 0, 0, 0, -1, 0, 1, 1, 0, 0, 0, 0, 0, 0, 0, 0, 1, -1, 0, 0, 0, 0, 0, 0, 0, 0, 0, 0, 0, 1, 0, 0, 0, 0, 0, -1, 1, 0, 0, 0, 0, 0, 0, -1, 0, 0, 0, 0, 0, 1, 1; 
0, 0, 0, 0, 0, 0, 0, 0, 0, 0, 0, 0, 0, 0, 0, 0, 0, 0, 0, 0, 0, -1, 0, 0, 0, 0, 0, 0, 0, 0, 0, 1, 0, 0, 0, 0, 0, 0, 0, 0, 0, 0, 0, -1, -1, 0, 0, 0, 0, 0, 0, 0, 0, 1, 0, 0, 0, 0, 0, 0, 0, 0, -1, 0, -1; 
0, 0, 0, 0, 0, 0, 0, 0, 0, 0, 0, 1, -1, 0, 0, 0, 0, 1, 0, 0, 0, 0, 0, 0, 0, -1, 0, 0, 0, 0, 0, 0, 0, 0, 0, 0, 0, 0, 0, 0, 0, 0, 0, 0, 1, 0, 0, 0, 0, 0, -1, 0, 0, 0, 0, 0, 0, 1, 0, 0, 0, 0, 0, 0, 0; 
0, 0, 0, 0, 0, 0, 0, 0, 0, 0, 0, 0, 0, 0, 0, -1, 1, 0, 0, 0, 0, 0, 0, 0, 0, 0, -1, 0, -1, 0, 0, 0, 0, 0, 0, 0, 0, 0, 0, 0, 0, 0, 0, 0, 0, 0, -1, 1, 0, 0, 0, 0, 0, 0, 0, 0, 0, 0, 1, 0, 1, 0, 0, 0, 1; 
0, 0, 0, 0, 0, 0, 0, 0, 0, 0, 0, 0, 0, 0, 0, 0, 0, 1, 1, 0, 0, 0, 0, 0, 0, 0, 0, 0, 0, 0, -1, 0, 0, 0, 0, 0, 0, 0, 0, 0, 0, 0, 0, 0, 0, 0, 0, 0, 0, -1, 0, 0, 0, 0, 0, 0, 0, 0, 0, 0, 0, 0, 1, 1, 0; 
0, 0, 0, 0, 0, 0, 0, 0, 0, 0, 0, 0, 0, 0, 0, 0, 0, -1, -1, 0, 0, 0, 0, 0, 0, 0, 0, 0, 0, 0, 1, 0, 0, 0, 0, 0, 0, 0, 0, 0, 0, 0, 0, 0, 0, 0, 0, 0, 0, 1, 0, 0, 0, 0, 0, 0, 0, 0, 0, 0, 0, 0, -1, -1, 0; 
0, 0, 0, 0, 0, 0, 0, 0, 0, 0, 0, 0, 0, 0, 0, 0, 0, 0, 0, -1, -1, 0, 0, 0, 0, -1, 0, 1, 0, -1, 0, 0, 0, 0, 0, 0, 0, 0, 0, 0, 0, 0, 0, 0, 0, 0, 0, 0, 0, 0, 0, 1, 1, 0, 0, 0, -1, 1, 0, 0, 0, 1, 0, 0, 2;  
0, 0, 0, 0, 0, 0, 0, 0, 0, 0, 0, 0, 0, 0, 0, 0, 0, 0, 0, 0, 0, 0, 0, -1, 1, 0, -1, 0, 0, 0, 0, 0, 0, 0, 0, 0, 0, 0, 0, 0, 0, 0, 0, 0, 0, 0, 0, 0, 0, 0, 0, 0, 0, 0, -1, 1, 0, 0, 1, -1, 0, 0, 0, 0, 0; 
0, 0, 0, 0, 0, 0, 0, 0, 0, 0, 0, 0, 0, 0, 0, 0, 0, 0, 0, 0, 0, 0, 0, 0, 0, 0, 0, -1, 1, 1, 0, 0, 0, 0, 0, 0, 0, 0, 0, 0, 0, 0, 0, 0, 0, 0, 0, 0, 0, 0, 0, 0, 0, 0, 0, 0, 0, 0, -1, 1, 0, 0, 0, 0, -1] \par}

Pari-Gp CALCULATIONS\\
$|b| = 2.890741236697218507543429035402903716418 E-27$\\
$L = 38.46960927036768465200292167581178343887$ \\
$|f'(a)^{-1}|_{\mathrm{sup}} = 8.212846275527759925085525656342053316915$ \\
$|f'(a)^{-1}|_{\mathrm{len}} = 10.32145710779244812406937753131330598443$ \\
$\frac{1}{2 {|f'(a)^{-1}|_{\mathrm{sup}}^{2}} L } = 0.0001926925132239904423664849871566682428236$\\
$\frac{1}{2 {|f'(a)^{-1}|_{\mathrm{len}}^{2}} L } = 0.0001220029142841818172845137711227723107218$ \\ 
\end{enumerate}
\subsection{Cusped Census}
We can apply the tests of Theorem~\ref{th} to every manifold in the SnapPea cusped census.  The results are found in the following theorem.
\begin{theorem}
Every manifold in the SnapPea cusped census has a complete hyperbolic structure.
\end{theorem}
A program was written in Perl~\cite{Perl} that issues commands to Snap to send tetrahedron shapes and filling equations for each manifold in the cusped census to an output file.  Then a Pari-Gp program reads the file, getting the needed data per manifold, and applies the template using this input.  The program then prints out the results.  The first run of this process determined that all but four manifolds, 5 168, 6 297, 7 1431 and 7 1927, have a complete hyperbolic structure.  The program rejected these four because each one, upon triangulation by Snap, had one tetrahedron shape parameter with an imaginary component that was effectively zero.  This was remedied by revising the original Perl program to process only these four manifolds, and including the ``randomize" command to get a different, acceptable triangulation.  The Pari-Gp program, also revised to process only these four manifolds, was then run using the second Perl output file.  The result was a determination that they also have a complete hyperbolic structure.

These programs can be adapted to give other information, such as the maximum value that $normb,$ the norm of $b,$ assumes over all the manifolds in the cusped census.  Call this $maxnormb$.  Similarly, for each manifold, we can ascertain the larger of the two values that $normb$ is compared to, and then the minimum of these maximum comparison values over all the manifolds in this census.  We do this because as long as $normb$ of a manifold is less than the larger of the two comparison values for that manifold, the manifold will have a complete hyperbolic structure.  Then if $normb$ of a manifold in the census is less than the smallest of these maximum comparison values over the whole census, that manifold is guaranteed to have a complete hyperbolic structure.  Call this minimum of maximum comparison values $minmaxvalue$.  It tells us the precision needed to evaluate a manifold in the census.  We have \\
$maxnormb$ = 1.717844093022015223183888589087321425164875899778 E-26 \\
$minmaxvalue$ = 0.00000147831677691814063380907736140260722549837777747014.\\ 
Thus, the approximate solution given by SnapPea, which is given to 10 digits but is computed to an internal precision of at least 15 significant digits, is sufficient for use as our $a_{1}, \ldots, a_{n}$.  It is interesting to see that the largest $normb$ is considerably smaller than the smallest comparison value over the entire cusped census.  The Perl programs and output files, as well as the Pari-Gp programs and log files, can be found at~\cite{HM}.  These Perl programs also include data with respect to a third test for a solution to the equations.  However, the third test yields smaller comparison values than the Kantorovich tests, so it has no effect on $minmaxvalue$.

\section*{Acknowledgements}
I want to thank Walter Neumann, my thesis adviser, for all his help and support since my first day at Columbia.  I am also very grateful to Joan Birman for her encouragement and advice over the years.  I am also appreciative of Chris Leininger's contribution to my examples.

\flushleft
\bibliography{thesisarticlebb}

\begin{thebibliography}{10}

\bibitem{Perl}
Comprehensive perl archive network.
\newblock available at http:/www.perl.com.

\bibitem{PGP}
Pari-gp.
\newblock available at http:/www.parigp-home.de.

\bibitem{B-P}
Riccardo Benedetti and Carlo Petronio.
\newblock {\em Lectures on Hyperbolic Geometry}.
\newblock Springer-Verlag, Berlin-Heidelberg-New York, 1992.

\bibitem{Choi}
Young-Eun Choi.
\newblock Positively oriented ideal triangulations on hyperbolic
  three-manifolds.
\newblock {\em Topology}, 43:1345--1371, 2004.

\bibitem{Edw}
Charles~Henry Edwards, Jr.
\newblock {\em Advanced Calculus of Several Variables}.
\newblock Dover Publcations,Inc., Mineola, N.Y., 1994.

\bibitem{Gab1}
David Gabai, Richard Meyerhoff, and Peter Milley.
\newblock Mom technology and volumes of hyperbolic 3-manifolds.
\newblock {\em arXiv.math.GT/0606072}, 2, 2006.

\bibitem{Gab2}
David Gabai, Richard Meyerhoff, and Peter Milley.
\newblock Minimum volume cusped hyperbolic three-manifolds.
\newblock {\em arXiv:0705.4325 [math.GT]}, 1, 2007.

\bibitem{Snap}
Oliver Goodman.
\newblock Snap.
\newblock available at http:/www.ms.unimelb.edu.au/~snap.

\bibitem{Kant}
John~H. Hubbard and Barbara~Burke Hubbard.
\newblock {\em Vector Calculus, Linear Algebra, and Differential Forms}.
\newblock Prentice Hall, Upper Saddle River, N.J., 1999.

\bibitem{Lein}
Christopher Leininger.
\newblock Small curvature surfaces in hyperbolic 3-manifolds.
\newblock {\em Journal of Knot Theory and It's Ramifications}, 15:379--411,
  2006.

\bibitem{HMdis}
Harriet Moser.
\newblock Doctoral dissertation, 2005.
\newblock available at http:/www.math.columbia.edu/~moser.

\bibitem{HM}
Harriet Moser.
\newblock Template, 2005.
\newblock available at http:/www.math.columbia.edu/~moser.

\bibitem{N}
Walter~D. Neumann.
\newblock Combinatorics of triangulations and the chern-simons invariant for
  hyperbolic 3-manifolds.
\newblock In {\em Topology 90}. Walter de Gruyter Co., Berlin-New York, 1992.

\bibitem{N-R}
Walter~D. Neumann and Alan~W. Reid.
\newblock Arithmetic of hyperbolic manifolds.
\newblock In {\em Topology 90}. Walter de Gruyter Co., Berlin-New York, 1992.

\bibitem{N-Z}
Walter~D. Neumann and Don Zagier.
\newblock Volumes of hyperbolic three manifolds.
\newblock {\em Topology}, 24(3):307--332, 1985.

\bibitem{Range}
Michael~R. Range.
\newblock {\em Holomorphic Functions and Integral Representations in Several
  Complex Variables}, volume 108 of {\em Graduate Texts in Mathematics}.
\newblock Springer-Verlag, Berlin-Heidelberg-New York, 1986.

\bibitem{Ratc}
John~G. Ratcliffe.
\newblock {\em Foundations of Hyperbolic Manifolds}, volume 149 of {\em
  Graduate Texts in Mathematics}.
\newblock Springer-Verlag, Berlin-Heidelberg-New York, 1948.

\bibitem{Thur2}
William~P. Thurston.
\newblock Hyperbolic structures on 3-manifolds i:deformation of acylindric
  manifolds.
\newblock {\em Ann. of Math.}
\newblock (preprint version).

\bibitem{Thur1}
William~P. Thurston.
\newblock The geometry and topology of 3-manifolds, 1979.
\newblock Princeton University mimeographed notes.

\bibitem{Thur}
William~P. Thurston.
\newblock {\em Three-Dimensional Geometry and Topology}.
\newblock Princeton University Press, Princeton, N.J., 1997.

\bibitem{Snapp}
Jeff Weeks.
\newblock Snappea.
\newblock available at http:/www.northnet.org/weeks.

\bibitem{Whit}
Hassler Whitney.
\newblock {\em Complex Analytic Varieties}.
\newblock Addison-Wesley Publishing Company, Reading,Mass.-Menlo
  Park,Ca.-London, 1972.

\end{thebibliography}
\bibliographystyle{plain}
\end{document}